\documentclass[a4paper,10pt]{amsart}
\usepackage{amsfonts}
\usepackage{anysize} \marginsize{1.3in}{1.3in}{1in}{1in}
 \usepackage{enumerate}
\usepackage{amsmath, amsthm}
\usepackage{amscd}
\usepackage{amssymb}
\usepackage{pdfsync}
\usepackage[title]{appendix}
\usepackage{xy}
\xyoption{all}
\usepackage{latexsym,graphicx}
\usepackage[latin1]{inputenc}
\usepackage{xspace}
\usepackage{array}
\usepackage{enumitem,kantlipsum}
\usepackage{newclude}
\usepackage[backref=page, bookmarksopen=true]{hyperref}
\hypersetup{pdfpagelabels,
                  plainpages=false,
                    colorlinks=true,       
                    linkcolor=blue}
\renewcommand*{\HyperDestNameFilter}[1]{\jobname-#1} 
\numberwithin{equation}{section}
\usepackage[a4paper, total={5in, 9in}, centering]{geometry}
\usepackage{todonotes}

\usepackage[nameinlink]{cleveref}

\newcommand{\noi}{\noindent}
\theoremstyle{plain}
\newtheorem{theor}{Theorem}[section]
\newtheorem{conj}[theor]{Conjecture}

\newtheorem{cor}[theor]{Corollary}

\theoremstyle{remark}
\newtheorem{rem}[theor]{Remark}

\newtheorem{Example}[theor]{Example}

\theoremstyle{plain}
\newtheorem{defi}[theor]{Definition}

\newcommand{\pos}{\textnormal{pos}}
\newcommand{\CC}{{\mathbb C}}
\newcommand{\RR}{{\mathbb R}}
\newcommand{\QQ}{{\mathbb Q}}
\newcommand{\ZZ}{{\mathbb Z}}
\newcommand{\VV}{{\mathbb V}}

\newcommand{\GG}{{\mathbf G}}
\newcommand{\HH}{{\mathbf H}}

\newcommand{\LL}{{\mathbf L}}

\newcommand{\PP}{{\mathbf P}}

\newcommand{\NN}{{\mathbb N}}

\newcommand{\Gr}{{\textnormal{Gr}}}
\newcommand{\HL}{\textnormal{HL}}
\newcommand{\ti}[1]{\mbox{$\tilde{#1} $}}
\newcommand{\wti}[1]{\mbox{$\widetilde{#1} $}}

\newcommand{\ol}{\overline}

\newcommand{\End}{\mathrm{End}}
\newcommand{\Hom}{\mathrm{Hom}}
\newcommand{\Res}{\mathrm{Res}}
\newcommand{\Sh}{{\rm Sh}}
\newcommand{\Gal}{{\rm Gal}}
\newcommand{\ad}{{\rm ad}}

\newcommand{\der}{{\rm der}}

\newcommand{\GL}{\textnormal{\textbf{GL}}}

\newcommand{\MT}{{\rm \bf MT}}

\newcommand{\alg}{\textnormal{alg}}
\newcommand{\proj}{{\mathbb P}}

\newcommand{\Sp}{\mathbf{Sp}}
\newcommand{\Aut}{\textnormal{Aut}}
\newcommand{\an}{\textnormal{an}}

\newcommand{\NL}{\textnormal{NL}}

\newcommand{\FS}{\mathfrak{S}}

\newcommand{\FH}{\mathfrak{H}}

\newcommand{\Fg}{\mathfrak{g}}

\newcommand{\cH}{{\mathcal H}}

\newcommand{\cE}{{\mathcal E}}

\newcommand{\cA}{{\mathcal A}}

\newcommand{\cS}{{\mathcal S}}
\newcommand{\cL}{{\mathcal L}}

\newcommand{\cX}{{\mathcal X}}
\newcommand{\cV}{{\mathcal V}}
\newcommand{\cO}{{\mathcal O}}
\newcommand{\cY}{{\mathcal Y}}

\newcommand{\bP}{\mathbf P}

\newcommand{\bN}{\mathbf N}

\newcommand{\bS}{\mathbf S}

\newcommand{\bT}{{\mathbf T}}

\newcommand{\oQ}{{\overline{\QQ}}}

\newcommand{\AAA}{{\mathbb A}}

\newcommand{\lto}{\longrightarrow}

\newcommand{\GSp}{\mathbf{GSp}}

\newcommand{\Lie}{\mathrm{Lie}}

\newcommand{\Imm}{\mathrm{Im}}

\newcommand{\Zar}{\textnormal{Zar}}

\newcommand{\Id}{\textnormal{Id}}

\newcommand{\Ad}{\textnormal{Ad}\,}

\newcommand{\prim}{\textnormal{prim}}
\newcommand{\dR}{\textnormal{dR}}
\newcommand{\hhom}{\textnormal{hom}}
\newcommand{\B}{\textnormal{B}}
\newcommand{\MHS}{\textnormal{MHS}}
\newcommand{\ul}{\underline}

\newcommand{\defin}{\textnormal{def}}
\newcommand{\Coh}{\textnormal{Coh}}
\newcommand{\codim}{\operatorname{codim}}
\newcommand{\atyp}{\textnormal{atyp}}
\newcommand{\typ}{\textnormal{typ}}
\newcommand{\fpos}{\textnormal{fpos}}

\newcommand{\BB}{\mathbf{B}}
\newcommand{\AH}{\textnormal{AH}}

\numberwithin{equation}{section}

\begin{document}

\title{Hodge theory, between algebraicity and transcendence}


\begin{abstract}
The Hodge theory of complex algebraic varieties is at heart a
transcendental comparison of two algebraic structures. 
We survey the recent advances bounding this transcendence, mainly due
to the introduction of o-minimal geometry as a natural framework for
Hodge theory. 
\end{abstract}

\maketitle


\section{Introduction}
Let $X$ be a smooth connected projective variety over $\CC$, and $X^\an$ its
associated compact complex manifold. Classical Hodge theory 
\cite{H51} states that the Betti (i.e. singular) cohomology 
group $H^k_\B(X^\an, \ZZ)$ is a \emph{polarizable $\ZZ$-Hodge structure of
  weight $k$}: there exists a canonical decomposition (called the
Hodge decomposition) of complex vector spaces
 \begin{equation*}
    H^k_{\B}(X^\an, \ZZ)\otimes_\ZZ \CC = \bigoplus_{p+q=k}
    H^{p,q}(X^\an)\;\; \textnormal{satisfying}\;\;
    \overline{H^{p,q}(X^\an)} = H^{q, p}(X^\an)
  \end{equation*}
and a $(-1)^k$-symmetric bilinear pairing 
  $
    q_k: H^k_{\B}(X^\an, \ZZ) \times
    H^k_{\B}(X^\an, \ZZ) \to \ZZ
  $ whose complexification makes the above decomposition orthogonal, and
  satisfies the
  positivity condition (the signs are complicated but are imposed to
  us by geometry)
  \begin{equation*}
    \mathrm{i}^{p-q} q_{k, \CC}(\alpha,
  \overline{\alpha})>0 \;\; \textnormal{for any nonzero} \; \alpha \in
  H^{p,q}(X^\an).
\end{equation*}
Deligne \cite{Del71} vastly generalized Hodge's result, showing that
the cohomology $H^k_{\B}(X^\an, \ZZ)$ of \emph{any} complex algebraic variety $X$ is functorially endowed with a slightly more 
general {\em mixed $\ZZ$-Hodge structure}, that makes, after tensoring with
$\QQ$, $H^k_{\B}(X^\an, \QQ)$
a successive extension of polarizable $\QQ$-Hodge structures, with weights between $0$ and $2k$. As mixed $\QQ$-Hodge structures form a
Tannakian category $\MHS_\QQ$, one can conveniently (although rather
abstractly) summarise the Hodge-Deligne
theory as functorially assigning to any complex algebraic variety
$X$ a $\QQ$-algebraic group: {\em the Mumford-Tate group $\MT_X$ of
  $X$}, defined as the Tannaka
group of the Tannakian subcategory $\langle H^\bullet_\B(X^\an, \QQ) \rangle$ of
$\MHS_\QQ$ generated by $H^\bullet_\B(X^\an, \QQ)$. The knowledge of
the group $\MT_X$ is equivalent to the knowledge of all {\em Hodge
  tensors} for the Hodge structure $H^\bullet_\B(X^\an, \QQ)$.

\medskip
These apparently rather innocuous semi-linear algebra statements are anything but
trivial. They have become the main tool for analysing the topology,
geometry and arithmetic of complex algebraic
varieties. Let's illustrate what we mean with regard to topology,
which we won't go into later. The existence of the Hodge
decomposition for smooth projective complex varieties, which holds more generally for compact K\"ahler
manifolds, imposes many constraints on the cohomology of such spaces, the most obvious
being that their odd Betti numbers have to be even. Such constraints
are not satisfied even by compact complex manifolds as simple
as the Hopf surfaces, quotients of $\CC^2 \setminus \{0\}$ by the
action of $\ZZ$ given by multiplication by $\lambda\not = 0$,
$|\lambda| \not =1$, whose first Betti number is
one. Characterising the homotopy types of compact K\"ahler manifolds
is an essentially open question, which we won't discuss here.

\medskip
The mystery of the Hodge-Deligne theory lies in the fact
that it is at heart \emph{not} an algebraic theory, but rather the
transcendental comparison of two algebraic structures. For simplicity
let $X$ be a smooth connected projective variety over $\CC$. The Betti cohomology $H^\bullet_\B(X^\an, \QQ)$ defines 
a $\QQ$-structure on the complex vector space of the algebraic de Rham cohomology
$H^\bullet_{\dR}(X/\CC):= H^\bullet (X, \Omega^\bullet_{X/\CC})$ under
the transcendental comparison isomorphism:
\begin{equation} \label{e1}
  \varpi: H^\bullet_{\dR}(X/\CC) \xrightarrow{\sim} H^\bullet(X^\an,
  \Omega^\bullet_{X^\an}) = \colon H^\bullet_{\dR}(X^\an, \CC) \xrightarrow{\sim} 
  H^\bullet_\B(X^\an, \QQ) \otimes_\QQ \CC\;\;,
\end{equation}
where the first canonical isomorphism is the comparison between algebraic and
analytic de Rham cohomology provided by GAGA, and the second one is provided by integrating complex
$\textnormal{C}^\infty$ differential forms over cycles (de Rham's
theorem). The Hodge filtration $F^p$ on $H^\bullet_\B(X^\an,
\QQ) \otimes_\QQ \CC$ is the image under $\varpi$ of the algebraic
filtration $F^p= \Imm(H^\bullet(X, \Omega^{\bullet \geq p}_{X/\CC})
\to H^\bullet_{\dR}(X/\CC))$ on the left hand side.

\medskip
The surprising power of the Hodge-Deligne theory lies in the fact that,
although the comparison between the two algebraic structures
is transcendental, this transcendence should be severely constrained, as
predicted for instance by the Hodge conjecture and the Grothendieck period
conjecture:

- For $X$ smooth projective, it is well-known that the cycle class
$[Z]$ of any codimension $k$ algebraic cycle on $X$ with $\QQ$
coefficients is a Hodge
class in the Hodge structure $H^{2k}(X^\an, \QQ)(k)$. Hodge \cite{H51}
famously conjectured that the converse holds true: any Hodge
class in $H^{2k}(X, \QQ)(k)$ should be such a cycle class.

- For $X$ smooth and defined over a number field $K \subset \CC$, its \emph{periods} are
the coefficients of the matrix of Grothendieck's
isomorphism (generalising  $(\ref{e1})$) 
\begin{equation*}
  \varpi: H^\bullet_{\dR}(X/K) \otimes_K \CC \xrightarrow{\sim}  H^\bullet_\B(X^\an, \QQ)
  \otimes_\QQ \CC
\end{equation*}
with respect to bases of
$H^\bullet_{\dR}(X/K)$ and $H^\bullet_\B(X^\an, \QQ)$. The
Grothendieck period conjecture (combined with the
Hodge conjecture) predicts that the transcendence degree of the field
$k_X \subset \CC$ generated by the periods of $X$ coincides with the
dimension of $\MT_X$.

\medskip
This tension between algebraicity and transcendence is perhaps best
revealed when considering Hodge theory {\em in families}, as developed
by Griffiths \cite{Grif}. Let $f: X \to S$ be a smooth projective morphism of smooth connected 
quasi-projective varieties over $\CC$. Its complex analytic fibers $X_s^\an$, $s \in
S^\an$, are diffeomorphic, hence their cohomologies 
  $\VV_{\ZZ, s}:= H^\bullet_\B(X_s^\an, \ZZ)$, $s \in S^\an$ are all isomorphic to a
  fixed abelian group $V_\ZZ$ and glue together into a locally constant
  sheaf $\VV_\ZZ:=R^\bullet {f^\an}_*\ZZ$ on $S^\an$. However the complex
  algebraic structure on $X_s$, hence also the Hodge structure on
  $\VV_{\ZZ, s}$, varies with $s$, making $R^\bullet f^\an_*\ZZ$ a variation of $\ZZ$-Hodge
  structures ($\ZZ$VHS) $\VV$ on $S^\an$, which can be naturally polarised. One easily checks that the
  Mumford-Tate group $\GG_s:= \MT_{{X_{s}}}$, $s \in
  S^\an$, is locally constant equal to the so-called \emph{generic Mumford-Tate
    group} $\GG$, outside of a meagre set $\HL(S, f)
  \subset S^\an$, the \emph{Hodge locus of the morphism $f$}, where it
  shrinks as exceptional Hodge tensors appear in
  $H^\bullet_\B(X_s^\an, \ZZ)$. The variation $\VV$ is completely
  described by its {\em period map}
  \begin{equation*}
    \Phi: S^\an \to \Gamma \backslash D.
  \end{equation*}
Here the period domain $D$ classifies all possible 
$\ZZ$-Hodge structure on the abelian group $V_\ZZ$, with a fixed
polarisation and Mumford-Tate
group contained in $\GG$; and $\Phi$ maps a point $s \in
S^\an$ to the point of $D$ parameterizing the polarized $\ZZ$-Hodge structure on $V_\ZZ$ defined by
$\VV_{\ZZ,s}$ (well-defined up to the action of the arithmetic group
$\Gamma:= G \cap
\GL(V_\ZZ)$).

\medskip
The transcendence of the comparison isomorphism $(\ref{e1})$ for each
fiber $X_s$ is embodied in the fact that the Hodge variety $\Gamma \backslash D$ is, in general, a mere complex
analytic variety not admitting any algebraic structure; and that the period map
$\Phi$ is a mere complex analytic map. On the other hand this transcendence
is sufficiently constrained so that the following corollary of the Hodge conjecture \cite{Weil} holds true, as proven by 
Cattani-Deligne-Kaplan \cite{CDK95}: the Hodge locus $\HL(S, f)$ is
a countable union of {\em algebraic} subvarieties of
$S$. Remarkably, their result is in fact valid for any polarized
$\ZZ$VHS $\VV$ on $S^\an$, not necessarily coming from geometry: the
Hodge locus $\HL(S, \VV^\otimes)$ is a countable union of algebraic
subvarieties of $S$.

\medskip
  In this paper we report on recent advances in the understanding of
  this interplay between algebraicity and transcendence in Hodge theory,
  our main object of interest being period maps $\Phi: S^\an \to
  \Gamma \backslash D$. The paper is written for non-experts: we
  present the mathematical objects involved, the questions, and the
  results but give only vague ideas of proofs, if any. It is organised as follows. After
  the \Cref{prelim} presenting the objects of Hodge
  theory (which the advanced reader will skip to refer to on occasion), we present in \Cref{tame} the main driving
  force behind the recent advances: although period maps are very rarely complex algebraic, their geometry is tame and does
  not suffer from any of the many possible pathologies of a 
  general holomorphic map. In model-theoretic terms, period maps are
  definable in the \emph{o-minimal structure} $\RR_{\an, \exp}$. In
  \Cref{transcendence}, we introduce the general format of
  \emph{bi-algebraic structures} for 
  comparing the algebraic structure on $S$ and the one on (the compact
  dual $\check{D}$ of) the period domain $D$. The heuristic provided
  by this format, combined with o-minimal geometry, leads to a powerful functional
  transcendence result: the Ax-Schanuel theorem for polarized
  $\ZZ$VHS. It also suggests to interprete variational Hodge theory as a special
  case of an {\em atypical intersection} problem. In \Cref{atypicalsection}
  we describe how this viewpoint leads to a stunning improvement of
  the result of Cattani, Deligne and Kaplan: in most cases $\HL(S,
  \VV^\otimes)$ is not only a countable union of algebraic varieties,
  but is actually algebraic on the nose (at least if we restrict to its
  components of positive period dimension). Finally in \Cref{arithm}
  we turn briefly to some arithmetic aspects of the theory. 

For the sake of simplicity we focus on the case of pure Hodge
structures, only mentioning the references dealing with the mixed
case. 

\subsection{Acknowledgments}
I would like to thank Gregorio Baldi, Benjamin Bakker, Yohan
Brunebarbe, Jeremy Daniel, Philippe
Eyssidieux, Ania Otwinowska, Carlos Simpson, Emmanuel Ullmo,
Claire Voisin, and Andrei Yafaev for many interesting discussions on
Hodge theory. I also thank Gregorio Baldi, Tobias Kreutz and Leonardo
Lerer for their comments on this text.

\section{Variations of Hodge structures and period maps} \label{prelim}

\subsection{Polarizable Hodge structures} \label{HS}
Let $n \in \ZZ$.
Let $R= \ZZ$, $\QQ$ or $\RR$. An \emph{$R$-Hodge structure $V$
  of weight $n$} is a finitely generated $R$-module $V_R$ together with one
of the following equivalent data: a bigrading $V_\CC(:= V_R \otimes_R
\CC)=\bigoplus_{p+q=n} V^{p,q}$, called the Hodge 
decomposition, such that $\overline{V^{p,q}} =V^{q,p}$ (the numbers
$(\dim V^{p,q})_{p+q=n}$ are called the Hodge numbers of $V$); or a decreasing filtration
 $F^\bullet$ of $V_\CC$, called the Hodge filtration, satisfying $F^p
 \oplus \ol{F^{n+1-p}}= V_\CC$. One goes from one to the other
    through $F^p = \bigoplus_{r\geq p} V^{r, n-r}$ and $V^{p,q} = F^p
    \cap \ol{F^q}$. The following group-theoretic description will be
    most useful to us: a Hodge structure is an $R$-module
    $V_R$ and a real algebraic representation $\varphi: 
    \bS \to \GL(V_\RR)$ whose restriction to $\GG_{m, \RR}$ is
    defined over $\QQ$. Here the Deligne torus $\bS$ denotes the real
    algebraic group $\CC^*$ of invertible matrices of the forms
    $\left( \begin{smallmatrix} a & -b 
  \\ b & a \end{smallmatrix} \right)$, which contains the diagonal
subgroup $\GG_{m, \RR}$. Being of weight $n$ is the requirement that
$\varphi_{|\GG_{m, \RR}}$ acts via the character $z \mapsto z^{-n}$. The space
$V^{p,q}$ is recovered as the eigenspace for the character $z\mapsto
z^{-p} {\ol{z}}^{-q}$ of $\bS(\RR)\simeq \CC^*$. A \emph{morphism of
Hodge structures} is a morphism of $R$-modules compatible with the
bigrading (equivalently, with the Hodge filtration or the
$\bS$-action).

\begin{Example}
 We write $R(n)$ for the unique $R$-Hodge structure of weight $-2n$,
called the Tate-Hodge structure of weight $-2n$, 
on the rank one free $R$-module $(2 \pi \mathrm{i})^n R \subset \CC$.
\end{Example}

One easily checks that the category of $R$-Hodge
structures is an abelian category (where the kernels and cokernels
coincide with the usual kernels and cokernels in the category of
$R$-modules, with the induced Hodge filtrations on their
complexifications), with natural tensor products $V \otimes W$ and
internal homs
$\hhom(V, W)$ (in particular duals $V^\vee: = \hhom(V, R(0))$). For
$R=\QQ$ or $\RR$ we obtain 
a Tannakian category, with an obvious exact
faithful $R$-linear tensor functor $\omega: (V_R, \varphi) \mapsto
V_R$. In particular $R(n)= R(1)^{\otimes n}$. If $V$ is an $R$-Hodge
structure we write $V(n):= V \otimes R(n)$ its $n$-th
Tate twist.

\medskip
If $V=(V_R, \varphi)$ is an $R$-Hodge structure of weight $n$, a 
\emph{polarisation for $V$} is a morphism of $R$-Hodge structures $q:
V^{\otimes 2} \lto R(-n)$ such that $(2 \pi \mathrm{i})^n q(x, \varphi(\mathrm{i})y)$
is a positive-definite bilinear form on $V_\RR$, called the
\emph{Hodge form} associated with the polarisation. If there exists a
polarisation for $V$ then $V$ is said \emph{polarizable}. One easily
checks that the category of polarizable $\QQ$-Hodge structures is
semi-simple.

\begin{Example} \label{ex1}
  Let $M$ be a compact complex manifold. If $M$
  admits a K\"ahler metric, the singular cohomology $H^n_{\B}(M, \ZZ)$ is
  naturally a $\ZZ$-Hodge structure of
  weight $n$, see \cite{H51}, \cite[Ch.6]{Voisin02}:
  \begin{equation*}
    H^n_{\B}(M, \ZZ)\otimes_\ZZ \CC = H^n_{dR}(M, \CC) = 
    \bigoplus_{p+q=n} H^{p,q} (M),
    \end{equation*}
where $H^\bullet_{dR}(M, \CC)$ denotes the de Rham cohomology
of the complex $(A^\bullet(M, \CC), d)$ of $\textnormal{C}^\infty$
differential forms on $M$, the first equality is the canonical isomorphism obtained by
integrating forms on cycles (de Rham theorem), and the complex vector
subspace $H^{p,q}(M)$ of $H^n_{\dR}(M, \CC)$ is generated by the $d$-closed forms of type
$(p,q)$, and thus satisfies automatically $\overline{H^{p,q}(M)} =
H^{q, p}(M)$. Although the second equality depends only on the complex
structure on $M$, its proof relies on the choice of a K\"ahler form
$\omega$ on $M$ through the following sequence of isomorphisms:
\begin{equation*}
  H^n_{dR}(M, \CC) \xrightarrow{\sim} \cH^n_{\Delta_{\omega}}(M) = \bigoplus_{p+q=n}
\cH_{\Delta_{\omega}}^{p,q} (M) \xrightarrow{\sim} 
\bigoplus_{p+q=n} H^{p,q} (M),
\end{equation*}
where $\cH^n_{\Delta_{\omega}} (M)$ denotes the vector space of $\Delta_\omega$-harmonic
differential forms on $M$ and $\cH_{\Delta_{\omega}}^{p,q} (M)$ its subspace of
$\Delta_\omega$-harmonic $(p,q)$-forms. The heart of Hodge theory is thus
reduced to the
statement that the Laplacian $\Delta_\omega$ of a K\"ahler metric
preserves the type of forms. The choice of a K\"ahler form $\omega$ on
$M$ also defines, through the hard
Lefschetz theorem \cite[Theor. 6.25]{Voisin02}, a polarisation {\em of
  the $\RR$-Hodge structure} $H^n(M, \RR)$, see
  \cite[Theor. 6.32]{Voisin02}. If $f: M\to N$ is any holomorphic map
  between compact complex manifolds admitting K\"ahler metrics
  then both $f^*: H^n_{\B}(N, \ZZ) \to H^n_{\B}(M,
  \ZZ)$ and the Gysin morphism
  $f_* : H^n_{\B}(M, \ZZ)  \to H^{n-2r}_{\B}(N,
  \ZZ)(-r)$ are morphism of $\ZZ$-Hodge structures, where $r= \dim M -
  \dim N$.
\end{Example}

\begin{Example} \label{ex2}
Suppose moreover that $M=X^\an$ is the compact complex manifold
analytification of a smooth projective variety $X$ over $\CC$. In that
case $H_{\B}^n(X, \ZZ)$ is a {\em polarizable $\ZZ$-Hodge
  structure}. Indeed, the K\"ahler class $[\omega]$ can be chosen as
the first Chern class of an ample line bundle on $X$, giving rise to a
rational Lefschetz decomposition and (after clearing denominators by
multiplying by a sufficiently large integer) to an integral
polarisation. Moreover the Hodge filtration $F^\bullet$ on $H^n_{\B}(X^\an,
\CC)$ can be defined algebraically: upon identifying $H^n_{\B}(X^\an,
\CC)$ with the algebraic de Rham cohomology $H^n_{\dR}(X/\CC):= H^n(X,
\Omega^\bullet_{X/\CC})$, the Hodge filtration is given by $F^p=
\Imm(H^n(X, \Omega^{\bullet \geq p}_{X/\CC}) \to H^n_{\B}(X^\an,
\CC))$. It follows that if $X$ is defined over a subfield $K$ of $\CC$,
then the Hodge filtration $F^\bullet$ on $H^n_{\B}(X^\an,
\CC) = H^n_{\dR}(X/K)\otimes_{K} \CC$ is defined over $K$.
\end{Example}

\begin{Example} \label{exAb}
 The functor which assigns to a complex abelian variety $A$ its
 $H^1_\B(A^\an, \ZZ)$ defines an equivalence of categories between abelian
 varieties and polarizable $\ZZ$-Hodge structures of weight $1$ and
 type $(1,0)$ and $(0,1)$.
 \end{Example}
\subsection{Hodge classes and Mumford-Tate group} \label{MumfordTate}

Let $R= \ZZ$ or $\QQ$ and let $V$ be an $R$-Hodge structure. A
\emph{Hodge class for $V$} is a vector in $V^{0,0} \cap V_\QQ =
F^0V_\CC \cap V_\QQ$. For instance, any morphism of $R$-Hodge
structures $f: V \to W$ defines a Hodge class in the internal $\hhom (V, W)$. Let $T^{m, n}V_\QQ$ denote the 
$\QQ$-Hodge structure $V_\QQ^{\otimes m} \otimes \hhom(V,
R(0))_\QQ^{\otimes n}$. A \emph{Hodge tensor for $V$} is a Hodge class
in some $T^{m, n}V_\QQ$. 

\medskip
The main invariant of an $R$-Hodge structure is its \emph{Mumford-Tate group}.
For any $R$-Hodge structure $V$ we denote by $\langle V
\rangle$ the Tannakian subcategory of the category of $\QQ$-Hodge structures
generated by $V_\QQ$ ; in
other words $\langle V \rangle$ is the smallest full subcategory
containing $V$, $\QQ(0)$ and
stable under $\oplus$, $\otimes$, and taking sub-quotients. If
$\omega_V$ denotes the restriction of the tensor functor $\omega$ to
$\langle V \rangle$, the functor $\Aut^\otimes (\omega_V)$ is representable by some closed
$\QQ$-algebraic subgroup $\GG_V \subset \GL(V_\QQ)$, called the Mumford-Tate
group of $V$, and $\omega_V$ defines
an equivalence of categories $\langle V \rangle \simeq
\textnormal{Rep}_\QQ \GG_V$. See \cite[II, 2.11]{Del82}.

\medskip
The Mumford-Tate group $\GG_V$ can also be characterised as the fixator in $\GL(V_\QQ)$ of the
Hodge tensors for $V$, or equivalently, writing $V= (V_R, \varphi)$, as the smallest $\QQ$-algebraic
subgroup of $\GL(V_\QQ)$ whose base change to $\RR$ contains the image
$\Imm \, \varphi$. In particular $\varphi$ factorises as $\varphi: \bS \to
\GG_{V, {\RR}}$. The group $\GG_V$ is thus connected, and reductive
if $V$ is polarizable. See \cite[Lemma2]{An92}.

\begin{Example}
  $\GG_{\ZZ(n)} = \GG_m$ if $n \not= 0$ and $\GG_{\ZZ(0)}= \{1\}$.
\end{Example}

\begin{Example}
Let $A$ be a complex abelian variety and let $V:= H^1_{\B}(A^\an, \ZZ)$ be
the associated $\ZZ$-Hodge structure of weight $1$. We write $\GG_A:=
\GG_V$. The choice of an
ample line bundle on $A$ defines a polarisation $q$ on $V$. On the one
hand, the endomorphism algebra $D:= \End^0(A) (:= \End(A) \otimes_\ZZ \QQ)$ is a
finite dimensional semisimple $\QQ$-algebra which, in view of
\Cref{exAb}, identifies with $\End(V_\QQ)^{\GG_{A}}$. Thus $\GG_A
\subset \GL_D(V_\QQ)$. On the other hand the polarisation $q$ defines
a Hodge class in $\hhom(V_\QQ^{\otimes 2}, \QQ(-1))$ thus $\GG_A$ has
to be contained in the group $\GSp(V_\QQ, q)$ of symplectic similitudes of
$V_\QQ$ with respect to the symplectic form $q$. Finally $\GG_A
\subset \GL_D(V_\QQ) \cap \GSp(V_\QQ, q)$.

If $A=E$ is an elliptic curve, it follows readily that either $D=
\QQ$ and $\GG_E=\GL_2$, or $D$ is an imaginary quadratic field
($E$ has complex multiplication) and $\GG_E= \bT_D$, the
$\QQ$-torus defined by $\bT_D(S) = (D\otimes_\QQ S)^*$ for any
$\QQ$-algebra $S$. 
\end{Example}

\subsection{Period domains and Hodge data} \label{period domains}
Let $V_\ZZ$ be a finitely generated abelian group $V_\ZZ$ of rank
$r$. Fix a positive integer $n$, a $(-1)^n$-symmetric
bilinear form $q_{\mathbf Z}$ on $V_\ZZ$ and a collection of non-negative
integers $(h^{p,q})$, ($p, q \geq 0$, $p+q=n$), such that $h^{p,q}= h^{q,p}$ and $\sum
h^{p,q}=r$. Associated with $(n, q_\ZZ, (h^{p,q}))$ we want to define a {\em period
  domain} $D$ classifying $\ZZ$-Hodge structures of weight $n$ on
$V_\ZZ$, polarized by $q_\ZZ$, and with Hodge numbers $h^{p,q}$. 
Setting $f^p = \sum_{r \geq p} h^{r, n-r}$ we first define the {\em compact
dual} $\check{D}$ parametrising the finite decreasing filtrations $F^\bullet$ on $V_\CC$
satisfying $(F^p)^{\perp_{q_{\ZZ}}}= F^{n+1-p}$ and $\dim F^p= f^p$. This
is a closed algebraic subvariety of the product of Grassmannians
$\prod_p \textnormal{\Gr}(f^p, V_\CC)$. The period domain $D\subset \check{D}^\an$ is the open subset where
the Hodge form is positive definite. If $\GG:=
\mathbf{GAut}(V_\QQ, q_\QQ)$ denotes the group of similitudes of $q_\QQ$, one
easily checks that $\GG(\CC)$ acts transitively on $\check{D}^\an$, which is
thus a flag variety for $\GG_\CC$; and that the connected
component $G:=\GG^\der(\RR)^+$ of the identity in the derived group $\GG^\der(\RR)$ acts transitively
on $D$, which identifies with an open $G$-orbit in 
$\check{D}$. If we fix a base point $o \in D$ and denote by $P$ and $M$
its stabiliser in $\GG(\CC)$ and $G$ respectively, the period
domain $D$ is thus the homogeneous space
\begin{equation*}
D= G/M \hookrightarrow \check{D}^\an = \GG(\CC)/P \;.
\end{equation*}
The group $P$ is a parabolic subgroup of $\GG(\CC)$. Its subgroup $M=
P \cap G$,
consisting of real elements, not only fixes the filtration
$F^\bullet_o$ but also the Hodge decomposition, hence the Hodge
form, at $o$. It is thus a compact subgroup of $G$ and $D$ is an open
elliptic orbit of $G$ in $\check{D}$.

\begin{Example}
 Let $n=1$, suppose that the only non-zero Hodge numbers are $h^{1,0}= h^{0,1}=g$,
 $q_\ZZ$ is a symplectic form and $D$ is the subset of
 $\textnormal{\Gr}(g, V_\CC)$ consisting of $q_\CC$-Lagrangian
 subspaces $F^1$ on which $\mathrm{i} q_\CC(u, \ol{u})$ is positive
 definite. In this case $\GG = \GSp_{2g}$, $G = \Sp_{2g}(\RR)$, $M= \mathbf{SO}_{2g}(\RR)$ is a
 maximal compact subgroup of the connected Lie group $G$ and $D= G/M$ is a
 bounded symmetric domain naturally biholomorphic to Siegel's upper
 half space $\mathfrak{H}_g$ of $g\times g$-complex symmetric
   matrices $Z= X+ \mathrm{i} Y$ with $Y$ positive
   definite. When $g=1$, $D$ is the Poincar\'e disk, biholomorphic to the
   Poincar\'e upper half space $\mathfrak{H}$.
 \end{Example}

More generally let $\GG$ be a connected reductive $\QQ$-algebraic group and let
$\varphi: \bS \to \GG_\RR$ be a real algebraic
morphism such that ${\varphi}_{|\GG_{m, \RR}}$ is defined over
$\QQ$. We assume that $\GG$ is the Mumford-Tate group of $\varphi$.
The \emph{period domain} (or
Hodge domain) $D$ associated with $\varphi: \bS \to \GG_\RR$
is the connected component of the $\GG(\RR)$-conjugacy class of
$\varphi: \bS \to \GG_\RR$ in $\Hom(\bS, \GG_\RR)$. Again, one easily checks that $D$ is an
open elliptic orbit of $G:=\GG^\der(\RR)^+$ in the compact dual flag variety $\check{D}^\an$, the
$\GG(\CC)$-conjugacy class of $\varphi_{\CC}\circ \mu:
\GG_{\textnormal{m}, \CC} \to \GG_\CC$, where $\mu:
\GG_{\textnormal{m}, \CC} \to \bS_\CC = \GG_{\textnormal{m}, \CC}
\times \GG_{\textnormal{m}, \CC}$ is the cocharacter $z \mapsto
(z,1)$. The pair $(\GG, D)$ is called a (connected) \emph{Hodge datum}. A
morphism of Hodge data $(\GG, D) \to (\GG', D')$ is a morphism $\rho:
\GG \to \GG'$ sending $D$ to $D'$. Any linear representation
$\lambda: \GG \to \GL(V_\QQ)$ defines a $\GG(\QQ)$-equivariant local
system $\check{\VV}_\lambda$ on $\check{D}^\an$. Moreover each point $x \in D$, seen as a
morphism $\varphi_x:\bS \to \GG_\RR$, defines a $\QQ$-Hodge structure
$V_x:= (V_\QQ, \lambda \circ \varphi_x)$. The
$\GG(\CC)$-equivariant filtration $F^\bullet \check{\cV}_\lambda:=
\GG^\ad(\CC) \times_{P, \lambda} F^\bullet V_{o,\CC}$ of the
holomorphic vector bundle $\check{\cV}_\lambda: = \GG^\ad(\CC)
\times_{P, \lambda} V_{o,\CC}$ on $\check{D}^\an$ induces the Hodge
filtration on $V_x$ for each $x \in D$. The Mumford-Tate group of
$V_x$ is $\GG$ precisely when $x \in D\setminus \bigcup \limits_{\tau}\tau(D')$, where $\tau$ ranges through the
countable set of morphisms of Hodge data $\tau: (\GG', D') \to (\GG,
D)$. The complex analytic subvarieties $\tau(D')$ of $D$ are called
the \emph{special subvarieties} of $D$.

\medskip
The following geometric feature of $\check{D}$ will be crucial for
us. The algebraic tangent bundle $T\check{D}$ naturally identifies, as
a $\GG_\CC$-equivariant bundle, with the quotient vector bundle $\check{\cV}_{\Ad}/ F^0
\check{\cV}_{\Ad}$, where $\Ad: \GG \to \GL(\Fg)$ is the adjoint
representation on the Lie algebra $\Fg$ of $\GG$. In particular it is
naturally filtered by the $F^i T\check{D}:= F^i \check{\cV}_{\Ad}/ F^0
\check{\cV}_{\Ad}$, $i \leq -1$. The subbundle $F^{-1}T\check{D}$ is
called the \emph{horizontal tangent bundle} of $\check{D}$.

\subsection{Hodge varieties}
Let $(\GG, D)$ be a Hodge datum as in \Cref{period domains}.
A {\em Hodge variety} is the quotient $\Gamma \backslash D$ of $D$ by
an arithmetic lattice $\Gamma$ of $\GG(\QQ)^+:= \GG(\QQ) \cap G$. It is thus naturally a
complex analytic variety, which is smooth if $\Gamma$ is
torsion-free. The {\em special subvarieties} of $\Gamma \backslash D$
are the images of the special subvarieties of $D$ under the projection
$\pi: D \to \Gamma \backslash D$ (one easily checks these are closed
complex analytic subvarieties of $\Gamma \backslash D$). For any algebraic
representation $\lambda: \GG \to \GL(V_\QQ)$, the
$\GG(\QQ)$-equivariant local system $\check{\VV}_\lambda$ as well
as the filtered holomorphic vector bundle $(\check{\cV}_\lambda,
F^\bullet)$ on $\check{D}$ are $G$-equivariant when restricted to $D$, hence descend to a
triple $(\VV_\lambda, (\cV_\lambda,F^\bullet), \nabla)$ on $\Gamma \backslash
D$. Similarly, the horizontal tangent bundle of $\check{D}$ defines
the \emph{horizontal tangent bundle} $T_h (\Gamma \backslash D)
\subset T(\Gamma \backslash D)$ of the Hodge variety $\Gamma
\backslash D$.

\subsection{Polarized $\ZZ$-variations of Hodge structures} \label{VHS}

Hodge theory as recalled in \Cref{HS} can be considered as the particular
case over a point of Hodge theory over an arbitrary base. Again, the motivation comes from
geometry. Let $f: Y \to B$ be a proper surjective complex analytic
submersion from a connected K\"ahler manifold $Y$ to a complex manifold $B$. It
defines a locally constant sheaf 
$\VV_\ZZ:= R^\bullet f_* \ZZ$ of finitely generated abelian groups on $B$, gathering the
cohomologies $H^\bullet_{\B}(Y_b, \ZZ)$, $b \in B$. Upon
choosing a base point $b_0 \in B$ the datum of $\VV_\ZZ$ is equivalent
to the datum of a {\em monodromy representation} $\rho: \pi_1(B, b_0)
\to \GL(\VV_{\ZZ, b_{0}})$. On the other hand, the de Rham incarnation
of the cohomology of the fibers of $f$ is the holomorphic
flat vector bundle $(\cV:= \VV_\ZZ \otimes_{\ZZ_{B}} \cO_B \simeq
R^\bullet f_* \Omega^\bullet_{Y/B}, \nabla)$,
where $\cO_B$ is the sheaf of holomorphic functions on $B$,
$\Omega^\bullet_{Y/B}$ is the relative holomorphic de Rham complex
and $\nabla$ is the Gauss-Manin connection. The Hodge filtration on each
$H^\bullet_{\B}(Y_b, \CC)$ is induced by the holomorphic subbundles
$F^p := R^\bullet f_* \Omega^{\bullet \geq p}_{Y/B}$ of $\cV$. The
Hodge filtration is usually not preserved by the connection, but
Griffiths \cite{Gri68} crucially observed that it satisfies the
\emph{transversality constraint} $\nabla F^p \subset \Omega^1_{B} \otimes_{\cO_{B}}  F^{p-1}
$.
More generally, \emph{a variation of $\ZZ$-Hodge
structures ($\ZZ$VHS)} on a connected complex manifold $(B, \cO_B)$ is a pair $\VV:= (\VV_\ZZ, 
F^\bullet)$, consisting of a locally constant sheaf of
finitely generated abelian groups $\VV_\ZZ$ on $B$ and a (decreasing)
filtration $F^\bullet$ of the holomorphic vector bundle $\cV:= \VV_\ZZ
\otimes_{\ZZ_{B}} \cO_B$ by holomorphic subbundles, called the Hodge filtration, satisfying the following
conditions: for each $b \in B$,
the pair $(\VV_b, F^\bullet_b)$ is a $\ZZ$-Hodge structure; and the
flat connection $\nabla$ on $\cV$ defined by $\VV_\CC$ satisfies Griffiths' transversality: 
\begin{equation} \label{GT}
\nabla F^\bullet \subset \Omega^1_{B} \otimes_{\cO_{B}} F^{\bullet
  -1}\;\;.
\end{equation}
A \emph{morphism $\VV \to \VV'$ of $\ZZ$VHSs} on $B$ is a morphism $f: \VV_\ZZ
\to \VV'_\ZZ$ of local systems such that the associated morphism of
vector bundles $f: \cV \to \cV'$ is compatible with 
the Hodge filtrations.
If $\VV$ has weight $k$, a \emph{polarisation of $\VV$} is a morphism $\mathrm{q}: \VV \otimes \VV
\to \ZZ_B(-k)$ inducing a polarisation on each $\ZZ$-Hodge structure
$\VV_b$, $b \in B$. In the geometric situation, such a polarisation exists if there exists an element $\eta\in H^2(Y, \ZZ)$ whose
restriction to each fiber $Y_b$ defines a K\"ahler class, for instance
if $f$ is the analytification of a smooth projective morphism of
smooth connected algebraic varieties over $\CC$.

\subsection{Generic Hodge datum and period map}
Let $S$ be a smooth connected quasi-projective variety over $\CC$ and let $\VV$ be a polarized $\ZZ$VHS on $S^\an$. Fix a base point $o \in S^\an$,
let $p: \wti{S^\an} \to S^\an$ be the corresponding universal cover and write
$V_\ZZ:= \VV_{\ZZ, o}$, $q_\ZZ:= \mathrm{q}_{\ZZ, o}$. The pulled-back
polarized $\ZZ$VHS $p^*\VV$ is canonically trivialised as $(\wti{S^\an} \times V_\ZZ,  (\wti{S^\an} \times V_\CC,
F^\bullet), \nabla = d, q_\ZZ)$. In 
\cite[7.5]{Del72}, Deligne proved that there exists a reductive $\QQ$-algebraic subgroup
$\iota: \GG
\hookrightarrow \GL(V_\QQ)$,
called the \emph{generic Mumford-Tate group} of $\VV$, such that, for
all points $\ti{s} \in \wti{S^\an}$, the Mumford-Tate group $\GG_{(V_\ZZ,
  F^\bullet_{\tilde{s}})}$ is contained in $\GG$, and is equal to
$\GG$ outside of a meagre set of $\wti{S^\an}$ (such points $\tilde{s}$
are said {\em Hodge generic} for $\VV$).  A closed irreducible
subvariety $Y \subset S$ is said Hodge generic for $\VV$ if it
contains a Hodge generic point. The setup of \Cref{period domains} is
thus in force. Without loss of generality we can assume
that the point $\tilde{o}$ is Hodge generic. Let $(\GG, D)$ be the Hodge
datum (called the \emph{generic Hodge datum} of $S^\an$ for $\VV$)
associated with the polarized Hodge structure $(V_\ZZ,
F^\bullet_{\tilde{o}})$. The $\ZZ$VHS $p^*\VV$ is completely described
by a holomorphic map $\tilde{\Phi}: \wti{S^\an} \to D$, which is 
naturally equivariant under the monodromy representation $\rho:
\pi_1(S^\an, o)\to \Gamma:= G \cap \GL(V_\ZZ)$, hence descends to a
holomorphic map $\Phi: S^\an \to \Gamma \backslash D$, called the
\emph{period map} of $S$ for $\VV$. We thus obtain the following
commutative diagram in the category of complex analytic spaces:
\begin{equation} \label{fundamental}
  \xymatrix{
    \wti{S^\an} \ar[r]^{\tilde{\Phi}} \ar[d]_{p} & D \; \ar[d]^{\pi}
    \ar@{^(->}[r] & \check{D}^\an\\
    S^\an \ar[r]^{\Phi} & \Gamma \backslash D\;. &
  }
\end{equation}

Notice that the pair $(\VV_\QQ, (\cV, F^\bullet))$ is the
pullback under $\Phi$ of the pair $(\VV_\iota, (\cV_\iota,
F^\bullet))$ on the Hodge variety $\Gamma \backslash D$ defined by the
inclusion $\iota: \GG \hookrightarrow \GL(V_\QQ)$. Griffiths' transversality
condition is equivalent to the statement that $\Phi$ is
\emph{horizontal}: $d\Phi(TS^\an) \subset T_h(\Gamma \backslash D)$. By
extension we call \emph{period map} any holomorphic, horizontal,
locally liftable map from $S^\an$ to a Hodge variety $\Gamma \backslash D$.

\medskip
The \emph{Hodge locus} $\HL(S, \VV^\otimes)$ of $S$ for $\VV$ is the subset of points $s \in 
S^{\an}$ for which the Mumford-Tate group $\GG_s$ is a strict subgroup
of $\GG$, or equivalently for which the Hodge structure $\VV_s$ admits more Hodge tensors
than the very general fiber $\VV_{s'}$. Thus 
\begin{equation} \label{HL}
\HL(S, \VV^\otimes) = \bigcup_{(\GG', D') \hookrightarrow (\GG, D)}
\Phi^{-1}(\Gamma'\backslash D') \;,
\end{equation}
where the union is over all strict Hodge subdata and
$\Gamma'\backslash D'$ is a slight abuse of notation for denoting the
projection of $D'\subset D$ to $\Gamma \backslash D$.

\medskip
Let $Y \subset S$ be a closed irreducible algebraic
subvariety $i: Y\hookrightarrow S$. Let $(\GG_Y, D_Y)$ be the generic
Hodge datum of the $\ZZ$VHS $\VV$ restricted to the smooth locus
of $Y$. The algebraic monodromy group $\HH_Y$ of $Y$ for
$\VV$ is the identity component of the Zariski-closure in $\GL(V_\QQ)$ of the monodromy of the
restriction to $Y$ of the local system $\VV_\ZZ$. It follows from
Deligne's (in the geometric case) and Schmid's (in general) ``Theorem of the fixed part'' and ``Semisimplicity
Theorem'' that $\HH_Y$ is a 
normal subgroup of the derived group $\GG_Y^{\textnormal{der}}$, see
\cite[Theorem 1]{An92}.

\section{Hodge theory and tame geometry} \label{tame}

\subsection{Variational Hodge theory between algebraicity and
  transcendence}

Let $S$ be a smooth connected quasi-projective variety over $\CC$ and
let $\VV=
(\VV_\ZZ, F^\bullet)$ be a
polarized $\ZZ$VHS on $S^\an$. Let $(\GG, D)$ be the generic Hodge
datum of $S$ for $\VV$ and let $\Phi: S^\an \to \Gamma \backslash D$
be the period map defined by $\VV$. 

The fact that Hodge theory is a transcendental theory
is reflected in the following facts:

\smallskip
- First, the triplets $(\VV_\lambda, (\cV_\lambda, F^\bullet),
\nabla)$ on $\Gamma \backslash D$ (for $\lambda: \GG \to \GL(V_\QQ)$ an algebraic
representation) do not in general satisfy Griffiths' transversality,
hence do not define a $\ZZ$VHS on $\Gamma \backslash D$. They do if
and only if $\VV$ is \emph{of Shimura type}, i.e $(\GG, D)$ is a (connected) {\em Shimura
datum} (meaning that the weight zero Hodge structures on the fibers of
$\VV_{\Ad}$ are of type $\{(-1,1), (0,0), (1, -1)\}$); or equivalently
if the horizontal tangent bundle $T_hD$ coincides with $TD$. In other
words: Hodge varieties are in general not classifying spaces for 
polarized $\ZZ$VHS.

\smallskip
- Second, and more importantly, the complex analytic Hodge variety
$\Gamma\backslash D$ is in general not algebraizable (i.e. it is not the
analytification of a complex quasi-projective variety). More precisely, let us write $D= G /M$
as in \Cref{period domains}. A classical property of elliptic orbits
like $D$ is that there exists a unique maximal compact
subgroup $K$ of $G$ containing $M$ \cite{gsc}. Supposing for simplicity that $G$ is a
real simple Lie group $G$, then $\Gamma \backslash D$ is algebraizable
only if $G/K$ is a hermitian symmetric domain and the projection $D
\to G/K$ is holomorphic or anti-holomorphic, see \cite{GRT}. 

\medskip
On the other hand this transcendence is severely constrained, as shown
by the following algebraicity results:

\smallskip
- If $(\GG, D)$ is \emph{of Shimura type}, then $\Gamma \backslash D=
\Sh^\an$ is the analytification of an algebraic variety, called a
Shimura variety $\Sh$ \cite{BB66}, \cite{De1},
\cite{De2}. In that case Borel \cite[Theor. 3.10]{Bor72} proved that the complex
analytic period map $\Phi: S^\an \to \Sh^\an$ is the analytification
of an algebraic map.

\smallskip
- Let $S \subset \ol{S}$ be a log-smooth compactification of $S$ by a
simple normal crossing divisor $Z$. Following Deligne \cite{Del}, the flat
holomorphic connection $\nabla$ on $\cV$ defines a canonical extension
$\ol{\cV}$ of $\cV$ to $\ol{S}$. Using
GAGA for $\ol{S}$, this defines an algebraic structure on $(\cV,
\nabla)$, for which the connection $\nabla$ is regular.
Around any point of $Z$, the complex manifold $S^\an$ is locally isomorphic to a product
$(\Delta^*)^k \times \Delta^l$ of punctured polydisks. Borel showed
that the monodromy representation $\rho: \pi_1(S^\an, s_o) 
\to \Gamma \subset \GG(\QQ)$ of $\VV$ is ``tame at infinity'', that
is, its restriction to $\ZZ^k = \pi_1((\Delta^*)^k \times \Delta^l)$ is
quasi-unipotent, see \cite[Lemma
(4.5)]{Schmid}. Using this result, Schmid showed that the Hodge
filtration $F^\bullet$ extends holomorphically to the Deligne
extension $\ol{\cV}$. This is the celebrated Nilpotent Orbit theorem \cite[(4.12)]{Schmid}. It
follows, as noticed by Griffiths \cite[(4.13)]{Schmid}), that the
Hodge filtration on $\cV$ comes from an algebraic filtration on the
underlying algebraic bundle, whether $\VV$ is of
geometric origin or not.

\smallskip
- More recently, an even stronger evidence came from
the study of Hodge loci. Cattani, Deligne and Kaplan proved the following
celebrated result (generalized to the mixed case in
\cite{BP1},  \cite{BP2}, \cite{BP3}, \cite{BPS}):

\begin{theor}[\cite{CDK95}] \label{CDK}
Let $S$ be a smooth connected quasi-projective variety over $\CC$ and
$\VV$ be a polarized $\ZZ$VHS over $S$. 
Then $\HL(S, \VV^\otimes)$ is a countable union of
closed irreducible algebraic subvarieties of $S$.
\end{theor}

In view of this tension between algebraicity and transcendence, it is
natural to ask if there is a framework, less strict than complex
algebraic geometry but more constraining than complex analytic
geometry, where to analyse period maps and explain its remarkable properties.

\subsection{O-minimal geometry}
Such a framework was in fact envisioned by Grothen\-dieck in \cite[\S
5]{Gro} under the name ``tame topology'', as a
way out of the pathologies of general topological spaces. Examples of
pathologies are Cantor sets, space-filling curves but also much
simpler objects like the graph $\Gamma:= \{(x, \sin \frac{1}{x}), \;
0<x\leq 1\} \subset \RR^2$: its closure $\ol{\Gamma}:= \Gamma \amalg
\textnormal{I}$, where $\textnormal{I}:= \{0\} \times [-1,1] \subset \RR^2$ is connected but not
arc-connected; $\dim (\ol{\Gamma} \setminus \Gamma) = \dim \Gamma$,
which prevents any reasonable stratification theory; and
$\Gamma \cap \RR$ is not ``of finite type''.  Tame geometry has been
developed by model theorists as o-minimal geometry, which studies
structures where every definable set has a finite geometric 
complexity. Its prototype is real semi-algebraic geometry, but it is
much richer. We refer to \cite{VDD} for a nice survey.

\begin{defi} A structure $\mathcal{S}$ expanding the real field is a collection $\mathcal{S}= (S_n)_{n \in \NN}$,
  where $S_n$ is a set of subsets of $\RR^n$ such that for every $n\in \NN$:
\begin{itemize}
\item[(1)] all algebraic subsets of $\RR^n$ are in $S_n$.
\item[(2)] $S_n$ is a boolean subalgebra of the power set of $\RR^n$
  (i.e. $S_n$ is stable by finite union, intersection, and complement).
\item[(3)] If $A\in S_n$ and $B \in S_m$ then $A \times B \in
  S_{n+m}$.
\item[(4)] Let $p: \RR^{n+1} \to \RR^n$ be a linear projection. If $A
  \in S_{n+1}$ then $p(A) \in S_n$.
\end{itemize}
The elements of $S_n$ are called the $\cS$-definable sets of
$\RR^n$. A map $f: A \to B$ between $\cS$-definable sets is said to be $\cS$-definable
if its graph is $\cS$-definable.
\end{defi}

A dual point of view starts from the functions, namely considers sets definable in a first-order
structure $\cS= \langle \RR, +, \times, <, (f_i)_{i \in I} \rangle$ where
$I$ is a set and the $f_i: \RR^{n_{i}} \to \RR$, $i \in I$, are functions. A subset $Z \subset
\RR^n$ is $\cS$-definable if it can be defined by a formula
$$ Z:= \{(x_1, \cdots , x_n) \in \RR^n \mid \phi(x_1, \cdots, x_n)
\; \textnormal{is true}\}\;\;,$$
where $\phi$ is a first-order formula that can be written using only
the quantifiers $\forall$ and $\exists$ applied to real variables; 
logical connectors; algebraic expressions written with the $f_i$;
the order symbol
$<$; and
fixed parameters $\lambda_i \in \RR$. When the set $I$ is empty the
$\cS$-definable subsets are the semi-algebraic sets. Semi-algebraic subsets
are thus always $\cS$-definable.

One easily checks that the composite of $\cS$-definable functions is
$\cS$-definable, as are the images and the preimages of $\cS$-definable
sets under $\cS$-definable maps. Using that the euclidean distance is
a real-algebraic function, one shows easily that the closure and
interior of an $\cS$-definable set are again $\cS$-definable.

\medskip
The following o-minimal axiom for a structure $\mathcal{S}$ guarantees the possibility of doing
geometry using $\cS$-definable sets as basic blocks.

\begin{defi}
A structure $\mathcal{S}$ is said to be o-minimal if $S_1$ consists
precisely of the finite unions of points and intervals
(i.e. the semi-algebraic subsets of $\RR)$. 
\end{defi}

\begin{Example}
The structure $\RR_{\sin}:= \langle \RR, +, \times, <, \sin \rangle$ is not
o-minimal. Indeed the infinite union of points $\pi \ZZ = \{x \in \RR
\mid \sin x=0\}$ is a definable subset of $\RR$ in this structure. 
\end{Example}

Any o-minimal structure $\cS$ has the following main tameness property:
given finitely many $\cS$-definable sets $U_1, \cdots, U_k \subset
\RR^n$, there exists a definable cylindrical cellular decomposition of
$\RR^n$ such that each $U_i$ is a finite union of cells. Such a
decomposition is defined inductively on $n$. For $n=1$ this is a
finite partition of $\RR$ into cells which are points or open
intervals. For $n>1$ it is obtained from a definable cylindrical cellular decomposition of
$\RR^{n-1}$ by fixing, for any cell $C\subset \RR^{n-1}$, finitely
many definable functions $f_{C, i}: C \to \RR$, $1\leq i \leq k_C$,
with $f_{C, 0}:= - \infty < f_{C,1 } < \cdots < f_{C, k_C} < f_{C, k_C
  +1}:= + \infty$, and defining the cells of $\RR^n$ as the graphs
$\{(x, f_{C, i}(x)), x \in C\}$, $1 \leq i \leq k_C$, and the bands $\{(x, f_{C, i}(x) < y < f_{C,
  i+1}(x)) , x \in C, y \in \RR$\}, $0 \leq i \leq k_C$, for all cells $C$ of
$\RR^{n-1}$.

\medskip
The simplest o-minimal structure is the structure $\RR_\alg$
consisting of semi-algebraic sets. It is too close to algebraic
geometry to be used for studying transcendence phenomena. 
Luckily much richer o-minimal geometries do
exist. A fundamental result of Wilkie, building on the
result of Khovanskii \cite{Khov} that any exponential set $\{(x_1,
\cdots, x_n) \in \RR^n\; |\; \;
P(x_1,  \dots, x_n, \exp(x_1), \dots, \exp(x_n))=0\}$ (where $P\in
\RR[X_1, \dots, X_n, Y_1, \dots, Y_n]$) has finitely many connected
components, states:
\begin{theor}[\cite{Wil}]
The structure $\RR_{\exp}:= \langle \RR, +, \times, <, \exp: \RR \to \RR
\rangle$ is o-minimal.
\end{theor}
In another direction, let us define
\begin{equation*} \RR_\an:=
\langle \RR,\, +, \,\times,\, <, \{f\} \; 
\textnormal{for} \, f \,\textnormal{restricted real analytic function}
\rangle,
\end{equation*}
where a function $f: \RR^n \to \RR$ is a restricted real analytic function if it is
zero outside $[0,1]^n$ and if there exists a real analytic function
$g$ on a neighbourhood of $[0,1]^n$ such that $f$ and $g$ are equal on
$[0,1]^n$. Gabrielov's result 
\cite{Gabrielov} that the difference of two subanalytic sets is
subanalytic implies rather easily that the structure $\RR_\an$ is o-minimal. The structure generated by two
o-minimal structures is not o-minimal in general, but Van den Dries
and Miller \cite{vdDM} proved that the structure $\RR_{\an, \exp}$
generated by $\RR_{\an}$ and $\RR_{\exp}$ is o-minimal.
This is the o-minimal structure which will be mainly used in the rest
of this text.

\medskip
Let us now globalize the notion of definable
set using charts:

\begin{defi}
  A definable topological space $\cX$ is the data of a Hausdorff
  topological space $\cX$, a \emph{finite} open covering $(U_i)_{1
    \leq i \leq k}$ of $\cX$, and homeomorphisms $\psi_i: U_i \to V_i \subset
  \RR^n$ such that all $V_i$, $V_{ij} := \psi_i(U_i \cap U_j)$ and
  $\psi_i \circ \psi_j^{-1}: V_{ij} \to V_{ji}$ are definable. As
  usual the pairs $(U_i , \psi_i)$ are called charts. A
  morphism of definable topological spaces is a continuous map which
  is definable when read in the charts. The definable site $\ul{\cX}$ of a definable topological
space $\cX$ has for objects definable open subsets $U \subset X$ and
admissible coverings are the finite ones. 
\end{defi}

\begin{Example} \label{definalg}
    Let $X$ be an algebraic variety over $\RR$. Then $X(\RR)$
    equipped with the euclidean topology carries a natural
    $\RR_\alg$-definable structure (up to isomorphism): one covers $X$
    by finitely many (Zariski) open affine subvarieties $X_i$ and take $U_i:= X_i(\RR)$
    which is naturally a semi-algebraic set. One easily check that any
    two finite open affine covers define isomorphic
    $\RR_{\alg}$-structures on $X(\RR)$.
If $X$ is an algebraic variety over $\CC$ then $X(\CC)=
    (\Res_{\CC/\RR} X)(\RR)$ carries thus a natural
    $\RR_\alg$-structure. We call this the {\em $\RR_\alg$-definabilization} of
    $X$ and denote it by $X^{\RR_\alg}$.
    \end{Example}

In the rest of this section, we fix an o-minimal structure $\cS$ and
write "definable" for $\cS$-definable. Given a complex algebraic
variety $X$ we write $X^\defin$ for the $\cS$-definabilization
$X^\cS$.

\subsection{O-minimal geometry and algebraization}
Why should an algebraic geometer care about o-minimal geometry?
Because o-minimal geometry provides strong algebraization results.

\subsubsection{Diophantine criterion}
The first algebraization result is the celebrated Pila-Wilkie theorem:
\begin{theor}[\cite{PW}] \label{PW}
  Let $Z\subset \RR^n$ be a definable set. We define $Z^\alg$ as the
  union of all connected {\em positive-dimensional} semi-algebraic subsets of
  $Z$. Then, denoting by $H: \QQ^n \to \RR$ the standard height function:
  \begin{equation*}
\forall \, \varepsilon>0, \quad \exists \, C_\varepsilon>0, \quad \forall
\, T>0, \quad \left| \{ x
  \in (Z\setminus Z^\alg) \cap \QQ^n, \; H(x) \leq T\} \right| <
C_\varepsilon T^\varepsilon.
\end{equation*}
\end{theor}

In words: if a definable set contains at least polynomially many rational points
(with respect to their height), then it contains a positive
dimensional semi-algebraic set! For instance, if $f: \RR \to \RR$
is a real analytic function such that its graph $\Gamma_f \cap [0,1]
\times [0,1]$ contains at least polynomially many rational points
(with respect to their height), then the function $f$ is real
algebraic \cite{BP}. This algebraization result is a
crucial ingredient in the proof of functional
transcendence results for period maps, see \Cref{transcendence}.

\subsubsection{Definable Chow and definable GAGA}
In another direction, algebraicity follows from the meeting of o-minimal geometry with complex
geometry. The motto is that o-minimal geometry is incompatible with the many pathologies
of complex analysis. As a simple illustration, let $f: \Delta^* \to \CC$ be
a holomorphic function, and assume that $f$ is definable (where we
identify $\CC$ with $\RR^2$ and $\Delta^*\subset 
\RR^2$ is semi-algebraic). Then $f$ does not have any essential
singularity at $0$ (i.e. $f$ is meromorphic). Otherwise, by the Big Picard
theorem, the boundary $\ol{\Gamma_f} \setminus \Gamma_f$ of its graph
would contain $\{0\} \times \CC$, hence would have the same real
dimension (two) as $\Gamma_f$, contradicting the fact that $\Gamma_f$
is definable.

\medskip
Let us first define a good notion of a definable topological space
``endowed with a complex analytic structure''. We identify $\CC^n$ with $\RR^{2n}$
by taking real and imaginary parts. Given $U \subset \CC^n$ a
definable open subset, let $\cO_{\CC^n}(U)$ denote the $\CC$-algebra
of holomorphic definable functions $U \to \CC$. The assignment $U
\rightsquigarrow \cO_{\CC^n}(U)$ defines a sheaf $\cO_{\CC^n}$ on $\ul{\CC^n}$ whose
stalks are local rings. Given a finitely generated ideal $I \subset
\cO_{\CC^n}(U)$, its zero locus $V(I) \subset U$ is definable and the
restriction $\cO_{V(I)}:= (\cO_U/I \cO_U)_{|\ul{V(I)}}$ define a sheaf
  of local rings on $\ul{V(I)}$.

  \begin{defi} \label{defiAnal}
    A definable complex analytic space is a pair $(\cX, \cO_{\cX})$
    consisting of a definable topological space $\cX$ and a sheaf 
    $\cO_{\cX}$ on $\ul{\cX}$ such that there exists a finite 
    covering of $\cX$ by definable open subsets $\cX_i$ on which
    $(\cX, \cO_{\cX})_{|\cX_{i}}$ is isomorphic to some $(V(I),
    \cO_{V(I)})$.
  \end{defi}

\cite[Theor. 2.16]{BBT} shows that
this is a reasonable definition: the sheaf $\cO_{\cX}$, in
analogy with the classical Oka's theorem, is a coherent sheaf of rings. Moreover one has a natural definabilization functor $(X,
  \cO_X) \rightsquigarrow (X^\defin, \cO_{X^\defin})$ from the category of
  separated schemes (or algebraic spaces) of finite type over $\CC$ to the category of
  definable complex analytic spaces, which induces a morphism $g:
  (\ul{X^\defin}, \cO_{X^\defin}) \to (\ul{X}, \cO_X)$ of locally
    ringed sites. 

\medskip
Let us now describe the promised algebraization results.
The classical Chow's theorem states that a closed complex
analytic subset $Z$ of $X^\an$ for $X$ smooth projective over $\CC$ is
in fact algebraic. This fails dramatically if $X$ is only
quasi-projective, as shown by the graph of the complex exponential in
$(\AAA^2)^\an$. However Peterzil and Starchenko, generalising \cite{FL}
in the $\RR_\alg$ case, have shown the following:

\begin{theor}[\cite{PS}, \cite{PS1}] \label{PS}
  Let $X$ be a complex quasi-projective variety and let $Z\subset X^\an$
  be a closed analytic subvariety. If $Z$ is definable in $X^\defin$ then $Z$
  is complex algebraic in $X$.
  \end{theor}

  Chow's theorem, which deals only with spaces, was extended to
  sheaves by Serre \cite{Serre}: when $X$ is proper, the analytification functor
  $(\cdot)^\an~:~\Coh(X)\to \Coh(X^\an)$ defines
  an equivalence of categories between the categories of coherent
  sheaves $\Coh(X)$ and $\Coh(X^\an)$. In the definable world, let $X$ be a separated scheme (or algebraic space) of finite type over
  $\CC$. Associating with a coherent sheaf $F$ on $X$ the coherent
  sheaf $F^\defin:= F \otimes_{g^{-1}\cO_X} \cO_{X^\defin}$ on the
  $\cS$-definabilization $X^\defin$ of $X$, one obtains a
  definabilization functor $(\cdot)^\defin :\Coh(X) \to
  \Coh(X^\defin)$. Similarly there is an analytification functor $\cX
  \rightsquigarrow \cX^\an$ from complex definable analytic spaces to
  complex analytic spaces, that induces a functor
  $(\cdot)^\an: \Coh(\cX) \to \Coh(\cX^\an)$. 

\begin{theor}[\cite{BBT}] \label{gaga}
For every separated algebraic space of finite type $X$, the
definabilization functor $(\cdot)^\defin: \Coh(X) \to \Coh(X^\defin)$ is
exact and fully faithful (but it is not necessarily essentially surjective). Its essential image is stable under
subobjects and subquotients.
\end{theor}

Using \Cref{gaga} and Artin's algebraization theorem for formal
modification \cite{Artin}, one obtains the following useful algebraization
result for definable images of algebraic spaces, which will be used in \Cref{Griffiths}:

\begin{theor}[\cite{BBT}] \label{images}
  Let $X$ be a separated
  algebraic space of finite type and let $\cE$ be a definable analytic
  space. Any proper definable analytic map $\Phi: X^\defin \to \cE$
  factors uniquely as $\iota \circ f^\defin$, where 
  $f: X\to Y$ is a proper morphism of separated algebraic spaces (of
  finite type) such that $\cO_Y \to f_* \cO_X$ is injective, and $\iota:
  Y^\defin \hookrightarrow \cE$ is a closed immersion of definable
  analytic spaces.
  \end{theor}

  \subsection{Definability of Hodge varieties}
Let us now describe the first result establishing that o-minimal geometry is potentially interesting for Hodge theory.

\begin{theor}[\cite{BKT}] \label{definPeriodDomain}
Any Hodge variety $\Gamma \backslash D$ can be naturally endowed with
a functorial structure $(\Gamma \backslash D)^{\RR_\alg}$ of $\RR_{\alg}$-definable
complex analytic space.
\end{theor}

Here ``functorial'' means that that any morphism $(\GG', D') \to (\GG, D)$ of Hodge data
induces a definable map $(\Gamma' \backslash D')^{\RR_\alg} \to (\Gamma
\backslash D)^{\RR_\alg}$ of Hodge varieties. Let us sketch the construction of
$(\Gamma \backslash D)^{\RR_\alg}$. Without loss of generality (replacing $\GG$
by its adjoint group if necessary) we can assume
that $\GG$ is semi-simple, $G =\GG(\RR)^+$. For simplicity let us
assume that the arithmetic lattice $\Gamma$ is torsion free. We choose a base point in $D =
G/M$. Notice that $G$ and $G/M \subset \check{D}^{\RR_{\alg}}$
are naturally endowed with a $G$-equivariant semi-algebraic structure,
making the projection $G \to G/M$ semi-algebraic. To define an $\RR_\alg$-structure on $\Gamma
\backslash (G/M)$, it is thus enough to find a semi-algebraic open 
fundamental set $F\subset G/M$ for the action of $\Gamma$ and to write $\Gamma
\backslash G/M= \Gamma \backslash F$, where the right hand side is the
quotient of $F$ by the closed \'etale semi-algebraic equivalence
relation induced by the action of $\Gamma$ on $D$. Here by fundamental
set we mean that the set of $\gamma \in \Gamma$ such that $\gamma F \cap F \not =
\emptyset$ is finite. We construct the fundamental set $F$ using the
reduction theory of arithmetic groups, namely
the theory of Siegel sets. Let $K$ be the unique maximal compact
subgroup of $G$ containing $M$. For any $\QQ$-parabolic $\bP$ of $\GG$
with unipotent radical $\bN$, the maximal compact subgroup $K$ of $G$
determines a real Levi $L\subset G$ which decomposes as $L=AQ$ where
$A$ is the center and $Q$ is semi-simple.  A semialgebraic Siegel set
of $G$ associated to $\bP$ and $K$ is then a set of the form
$\FS=U(aA_{>0})W$ where $U\subset \bN(\RR)$, $W\subset QK$ are
bounded semialgebraic subsets, $a\in A$, and $A_{>0}$ is the cone
corresponding to the positive root chamber. By a Siegel set of $G$
associated to $K$ we mean a semialgebraic Siegel set associated
to $\bP$ and $K$ for some $\QQ$-parabolic $\bP$ of $\GG$. Suppose now that $\Gamma \subset G$
is an arithmetic group. A fundamental result of Borel
\cite{bor} states that there exists finitely many Siegel
sets $\FS_i \subset G$, $1\leq i \leq s$, associated with
$K$, whose images in $\Gamma \backslash G/K$ cover the
whole space; and such that for any $1\leq i\not =j \leq s$, the set of
$\gamma \in \Gamma$ such that $\gamma \FS_i \cap \FS_j \not =
\emptyset$ is finite. We call the images $\FS_{i, D}:= \FS_i/M$ {\em Siegel sets for
  $D$}. Noticing that these Siegel sets for $D$ are semi-algebraic in
  $D$, we can take $F= \coprod_{i=1}^s \FS_{i, D}
$. It is not difficult to show that the $\RR_\alg$-structure thus
constructed is independent of the choice of the base point $eM \in
G/M$. The functoriality follows from a non-trivial property of Siegel
sets with respect to morphisms of algebraic groups, due to Orr \cite{Orr}.

\subsection{Definability of period maps}
Once \Cref{definPeriodDomain} is in place, the following result shows
that o-minimal geometry is a natural framework for Hodge theory:

\begin{theor}[\cite{BKT}] \label{definPeriodMap}
Let $S$ be a smooth connected complex quasi-projective variety. Any period map $\Phi:
S^\an \to \Gamma \backslash D$ is the analytification of a morphism $\Phi: S^{\RR_{\an, \exp}} \to (\Gamma
\backslash D)^{\RR_{\an, \exp}}$ of $\RR_{\an, \exp}$-definable complex
analytic spaces, where the
$\RR_{\an,\exp}$-structures on $S(\CC)$ and $\Gamma \backslash D$
extend their natural $\RR_\alg$-structures defined in \Cref{definalg}
and \Cref{definPeriodDomain} respectively.
\end{theor}

In down to earth terms, this means that we can cover $S$ by finitely
many open affine charts $S_i$ such that $\Phi$ restricted to
$(\Res_{\CC/\RR} S_i) (\RR)= S_i(\CC)$ and read in a chart of
$\Gamma \backslash D$ defined by a Siegel set of $D$, can be written using only
real polynomials, the real exponential function, and restricted real
analytic functions! This statement is already non-trivial when $S= \Sh$
is a Shimura variety and $\Phi^\an : S^\an \to \Gamma \backslash D$ is
the identity map coming from the uniformization $\pi: D \to S^\an$ of
$S^\an$ by the hermitian symmetric domain $D= G/K$. In that case the
$\RR_\alg$-definable varieties $\Sh^{\RR_{\alg}}$ and $(\Gamma
\backslash D)^{\RR_{\alg}}$ are not isomorphic but
\Cref{definPeriodMap} claims that their $\RR_{\an,\exp}$-extensions
$\Sh^{\RR_{\an, \exp}}$ and $(\Gamma
\backslash D)^{\RR_{\an, \exp}}$ are. 
This is equivalent to showing that the restriction $\pi_{|\FS_D}: \FS_D \to
S^{\RR_{\an, \exp}}$ to a Siegel set for $D$ can be written using only
real polynomials, the real exponential function, and restricted real
analytic functions. This is a nice exercise on the $j$-function when
$\Sh$ is a modular curve, was
done in \cite{PS2} and \cite{PT} for $\Sh= \cA_g$, and \cite{KUY} in general. 

\medskip
Let us sketch the proof of \Cref{definPeriodMap}. We choose a
log-smooth compactification of $S$, hence providing us with a
definable cover of $S^{\RR_{\an}}$ by punctured polydisks
$(\Delta^*)^k \times \Delta^l$. We are reduced to showing that
the restriction of $\Phi$ to such a punctured polydisk is $\RR_{\an,
  \exp}$-definable. This is clear if $k=0$, as in this case $\varphi:
\Delta^{k+l} \to \Gamma \backslash D$ is even
$\RR_{\an}$-definable. For $k>0$, let $\textnormal{e}:\exp(2\pi \mathrm{i}
\cdot): \FH \to \Delta^*$ be 
the universal covering map. Its restriction to a sufficiently large
bounded vertical strip $V:= [a,b] \times ]0, +\infty[ \subset \FH=
\{x+\mathrm{i} y, \; y>0\}$ is $\RR_{\an, \exp}$-definable. Considering
the following commutative diagram:
\begin{equation*}
  \xymatrix{
V^k \times \Delta^l \ar[rr]^{\tilde{\Phi}} \ar[d]_{\textnormal{e}} &  &D \ar[d]^\pi&
  F \ar@{_(->}[l]\\
  (\Delta^*)^k \times \Delta^l \, \ar@{^(->}[r] & S^{\an}
  \ar[r]_\Phi & \Gamma \backslash D \;,&
}
\end{equation*}
it is thus enough to show that $ \pi \circ \tilde{\Phi}: V^k \times \Delta^l  \to
\Gamma \backslash D$ is $\RR_{\an, \exp}$-definable.

Let the coordinates of $(\Delta^*)^k \times \Delta^l$ be $t_i$,
$1\leq i \leq k+l$,
those of $\FH^k$ be $z_i$, $1\leq i \leq k$, so that $\mathrm{e}(z_i)= t_i$. Let $T_i$ be
the monodromy at infinity of $\Phi$ around the hyperplane $(z_i=0)$,
boundary component of $\ol{S}\setminus S$. By Borel's theorem $T_i$ is
quasi-unipotent. Replacing $S$ by a finite \'etale cover we can
without loss of generality assume that each $T_i= \exp(N_i)$, with $N_i \in
\Fg$ nilpotent. The Nilpotent Orbit Theorem of Schmid is
equivalent to saying that $\tilde{\Phi}: V^k \times \Delta^l \to D$
can be written as $\tilde{\Phi}(z_1, \dots, z_k, t_{k+1}, \dots ,
t_{k+l}) = \exp(\sum_{i=1}^k z_i N_i) \cdot \Psi(t_1, \dots,
t_{k+l})$ for $\Psi: \Delta^k \times \Delta^l \to \check{D}^\an$ a
holomorphic map. On the one hand $\Psi$ is $\RR_{\an}$-definable as a
function of the variables $t_i$, hence $\RR_{\an, \exp}$-definable as a
function of the variables $z_i$, $1\leq i \leq k$, and the variables
$t_j$, $k+1\leq j \leq k+l$ . On the other hand $\exp(\sum_{i=1}^k z_i N_i)
\in \GG(\CC)$ is polynomial in the variables $z_i$, as the monodromies
$N_i$ are nilpotent and commute pairwise. As the action of $\GG(\CC)$
on $\check{D}$ is algebraic, it follows that $\tilde{\Phi}: V^k \times \Delta^l  \to
D$ is $\RR_{\an, \exp}$-definable. The proof of \Cref{definPeriodMap}
is thus reduced to the following, proven by Schmid when $k=1, l=0$ \cite[5.29]{Schmid}:
\begin{theor}[\cite{BKT}]
The image $\tilde{\Phi}(V^k \times \Delta^l)$ lies in a finite union
of Siegel sets of $D$.
\end{theor}

\noindent
This can be interpreted as showing that, possibly after passing to a
definable cover of $V^k \times \Delta^l$, the Hodge form of
$\tilde{\Phi}$ is Minkowski reduced with respect to a flat frame. This
is done using the hard analytic theory of Hodge forms estimates for
degenerations of variations of Hodge structure, as in
\cite[Theor. 3.4.1 and 3.4.2]{Ka} and \cite[Theor. 5.21]{CKS}. 

\begin{rem}
\Cref{definPeriodDomain} and \Cref{definPeriodMap} have been extended
to the mixed case in \cite{BBKT}.
\end{rem}

\subsection{Applications}

\subsubsection{About the Cattani-Deligne-Kaplan theorem}
As a corollary of \Cref{definPeriodMap} and \Cref{PS} one obtains
  the following, which, in view of (\ref{HL}), implies immediately
  \Cref{CDK}:
  
  \begin{theor}[\cite{BKT}] \label{algebraicity}
    Let $S$ be a smooth quasi-projective complex variety. Let $\VV$
    be a polarized $\ZZ$VHS on $S^\an$ with period map
    $\Phi: S^\an \to \Gamma \backslash D$. For any
    special subvariety $\Gamma' \backslash D' \subset \Gamma
    \backslash D$ its preimage $\Phi^{-1}(\Gamma' \backslash D')$ is a
    finite union of irreducible algebraic subvarieties of $S$.
  \end{theor}

  Indeed it follows from \Cref{definPeriodDomain} that $\Gamma'
  \backslash D'$ is definable in $(\Gamma
  \backslash D)^{\RR_{\alg}}$. By \Cref{definPeriodMap} its preimage
  $\Phi^{-1}(\Gamma' \backslash D')$ is definable
  in $S^{\RR_{\an, \exp}}$. As $\Phi$ is holomorphic and $\Gamma' \backslash D' \subset \Gamma
    \backslash D$ is a closed complex analytic subvariety,
    $\Phi^{-1}(\Gamma' \backslash D')$ is also a closed complex
    analytic subvariety of $S^\an$. By \Cref{PS} it is thus algebraic
    in $S$.

    \begin{rem}
      \Cref{algebraicity} has been extended to the mixed case in
      \cite{BBKT}, thus recovering \cite{BP1}, \cite{BP2}, \cite{BP3},
      and \cite{BPS}.
    \end{rem}

Let $Y \subset S$ be a closed irreducible algebraic subvariety. Let $(\GG_Y, D_Y) \subset (\GG, D)$ be the generic Hodge datum of $\VV$ restricted to the
smooth locus of of $Y$. There exist a smallest Hodge subvariety
$\Gamma_Y\backslash D_Y$ of $\Gamma \backslash D$ containing
$\Phi(Y^\an)$. The following terminology will be convenient:

 \begin{defi} \label{special}
      Let $S$ be a smooth quasi-projective complex variety. Let $\VV$
    be a polarized $\ZZ$VHS on $S^\an$ with period map
    $\Phi: S^\an \to \Gamma \backslash D$. A closed irreducible
    subvariety $Y \subset S$ is called a \emph{special subvariety} of $S$
    for $\VV$ if it coincides with an irreducible component of the preimage
    $\Phi^{-1}(\Gamma_Y\backslash D_Y)$.

   Equivalently, a special subvariety of $S$ for $\VV$ is a closed irreducible
algebraic subvariety $Y \subset S$ maximal among the closed
irreducible algebraic subvarieties $Z$ of $S$ such that the generic
Mumford-Tate group $\GG_Z$ of $\VV_{|Z}$ equals $\GG_Y$.
      \end{defi}

 \subsubsection{A conjecture of Griffiths} \label{Griffiths}
Combining \Cref{definPeriodMap} this time with 
\Cref{images} leads to a proof of 
an old conjecture of Griffiths \cite{Grifb}, claiming that the image
of any period map has a natural structure of quasi-projective variety
(Griffiths proved it when the target Hodge variety is compact):

\begin{theor}[\cite{BBT}]
  Let $S$ be a smooth connected quasi-projective complex variety and
  let $\Phi: S^\an \to \Gamma \backslash D$ be a period map. There
  exists a unique dominant morphism of complex algebraic varieties $f:
  S\to T$, with $T$ quasi-projective, and a closed complex analytic immersion $\iota: T^\an \hookrightarrow \Gamma
  \backslash D$ such that $\Phi = \iota \circ f^\an$.
\end{theor}

Let us sketch the proof.
As before, let $S\subset \ol{S}$ be a
log-smooth compactification by a simple normal crossing divisor
$Z$. It follows from a result of Griffiths \cite[Prop.9.11i)]{Grif}
that $\Phi$ extends to a {\em proper} period map over the components
of $Z$ around which the monodromy is finite. Hence, without loss of
generality, we can assume that $\Phi$ is proper. 
The existence of $f$ in the category of algebraic spaces
then follows immediately from \Cref{definPeriodMap} and
\Cref{images} (for $\cS = \RR_{\an, \exp}$). The proof that $T$ is in fact quasi-projective exploits
a crucial observation of Griffiths that $\Gamma\backslash D$ carries a
positively curved $\QQ$-line bundle $\cL:= \bigotimes_p \det(F^p)$. This
line bundle is naturally definable on $(\Gamma \backslash D)^\defin$. Using
the definable GAGA \Cref{gaga}, one shows that its restriction to
$T^\defin$ comes from an algebraic $\QQ$-line bundle $L_T$ on $T$,
which one manages to show to be ample.

\section{Functional transcendence} \label{transcendence}

\subsection{Bi-algebraic geometry}
As we saw, Hodge theory, which compares the Hodge filtration
on $H^\bullet_{\dR}(X/\CC)$ with the rational structure on
$H^\bullet_{\B}(X^\an, \CC)$, gives rise to variational Hodge theory,
whose fundamental diagram~(\ref{fundamental}) compares the
algebraic structure of $S$ with the algebraic structure on the dual
period domain $\check{D}$. As such, it is a partial answer to one of
the most classical problem of complex algebraic geometry: the
transcendental nature of the topological universal cover of complex
algebraic varieties. If $S$ is a connected complex algebraic variety,
the universal cover $\widetilde{S^\an}$ has usually no algebraic
structure as soon as the topological fundamental group $\pi_1(S^\an)$
is infinite. As an aside, let us mention an interesting conjecture of
K\'ollar and Pardon \cite{KP12}, predicting that if $X$ is a normal
projective irreducible complex variety whose universal cover
$\wti{X^\an}$ is biholomorphic to a semialgebraic open subset of an
algebraic variety then $\wti{X^\an}$ is biholomorphic to $\CC^n \times
D \times F^\an$, where $D$ is a bounded symmetric domain and $F$ is a
normal, projective, irreducible, topologically simply connected,
complex algebraic variety. We want to think of
variational Hodge theory as an attempt to provide a partial
\emph{algebraic uniformization}: the period map emulates an algebraic
structure on $\widetilde{S^\an}$, modeled on the flag variety
$\check{D}$. The remaining task is then to describe the transcendence properties of
the complex analytic uniformization map $p: \wti{S^\an} \to
S^\an$ with respects to the emulated algebraic structure on
$\wti{S^\an}$ and the algebraic structure $S$ on $S^\an$. A few years
ago, the author \cite{K17}, together with Ullmo and Yafaev
\cite{KUY2}, introduced a convenient format for studying such
questions, which encompasses many classical transcendence problems and
provides a powerful 
heuristic. 

\begin{defi}
A bi-algebraic structure on a connected quasi-projective variety
$S$ over $\CC$
is a pair 
\begin{equation*}
  (f: \wti{S^\an} \to Z^\an, \quad \rho: \pi_1(S^\an) \to \Aut(Z))
\end{equation*}
where $Z$ denotes an algebraic variety (called the {\em algebraic
  model} of $\wti{S^\an}$),  $\Aut(Z)$ is its group of algebraic
automorphisms, $\rho$ is a group morphism
(called the monodromy representation) and 
$f$ is a $\rho$-equivariant holomorphic map (called the developing
map).
\end{defi}

An irreducible analytic subvariety $Y \subset \wti{S^\an}$ is said to be an
{\em algebraic subvariety of $\wti{S^\an}$} for the bi-algebraic
structure $(f, \rho)$ if $Y$ is an analytic
irreducible component of $f^{-1}(\overline{f(Y)}^\Zar)$ (where
$\overline{f(Y)}^\Zar$ denotes the Zariski-closure of $f(Y)$ in $Z$).
An irreducible algebraic subvariety $Y\subset \wti{S^\an}$, resp. $W
\subset S$, is said to be {\em bi-algebraic} if $p(Y)$ is an
algebraic subvariety of $S$, resp. any (equivalently one) analytic irreducible component of $p^{-1}(W)$ is an
irreducible algebraic subvariety of $\wti{S^\an}$.
The bi-algebraic subvarieties of $S$ are precisely the ones where the
emulated algebraic structure on $\wti{S^\an}$ and the one on $S$
interact non-trivially.

\begin{Example} \label{biExample}

  \smallskip
  \noindent
 (a) {\it tori:} $S= (\CC^*)^n$. The uniformization map is the multi-exponential
$$p:= (\exp(2\pi i \cdot), \dots,\exp(2\pi i \cdot)) : \CC^n \to
(\CC^*)^n,$$ and $f$ is the identity morphism of $\CC^n$.
An irreducible algebraic subvariety $Y \subset \CC^n$ (resp. $W \subset (\CC^*)^n$) is bi-algebraic if 
and only if $Y$ is a translate of a rational linear subspace of $\CC^n= \QQ^n
\otimes_\QQ \CC$ (resp. $W$ is a translate of a subtorus of
$(\CC^*)^n$). 

\smallskip
\noindent
(b) {\it abelian varieties:} $S=A$ is a complex abelian variety of dimension $n$.
Let $p: \Lie \,A \simeq \CC^n \to A$ be the uniformizing map of a complex abelian
variety $A$ of dimension $n$. Once more $\wti{S^\an} =  \CC^n$ and $f$ is the identity morphism. 
One checks easily that an irreducible algebraic subvariety $W \subset A$ is bi-algebraic if
and only if $W$ is the translate of an abelian subvariety of $A$.

\smallskip
\noindent
(c) {\it Shimura varieties:} Let $(\GG, D)$ be a
Shimura datum. The quotient $S^\an = \Gamma \backslash D$ (for $\Gamma \subset G:=
\GG^\der(\RR)^+$ a congruence torsion-free lattice) is the complex
analytification of a (connected) Shimura variety $\Sh$, defined over a
number field (a finite extension of the reflex field of
$(\GG, D)$). And $f$ is the open embedding $D \hookrightarrow
\check{D}^\an$.
\end{Example}

Let us come back to the case of the bi-algebraic structure on $S$
\begin{equation*}
  (\tilde{\Phi}: \wti{S^\an}\to \check{D}^\an, \; \rho: \pi_1(S^\an) \to \Gamma \subset
\GG(\QQ))
\end{equation*} defined by a polarized $\ZZ$VHS $\VV$ and its period map $\Phi:
S^\an \to \Gamma \backslash D$ with monodromy $\rho: \pi_1(S^\an) \to
\Gamma \subset \GG(\QQ)$ (in fact all the examples above are of this form if we
consider more generally graded-polarized variations of mixed
$\ZZ$-Hodge structures). What are its bi-algebraic subvarieties? To
answer this question, we need to define the \emph{weakly special}
subvarieties of $\Gamma \backslash D$, as either a special subvariety or a subvariety of the form $$\Gamma_\HH \backslash
D_H \times \{t\} \subset  \Gamma_\HH \backslash
D_H\times \Gamma_\LL \backslash
D_L \subset \Gamma \backslash D\;\;,$$ where $(\HH \times \LL, D_H
\times D_L)$ is a Hodge subdatum of $(\GG^\ad, D)$ and $\{t\}$ is a Hodge
generic point in $\Gamma_\LL \backslash
D_L$. 
Generalising \Cref{algebraicity}, the preimage under $\Phi$
of any weakly special subvariety of $\Gamma \backslash D$ is an algebraic
subvariety of $S$, \cite{KO}. An irreducible component of such a preimage is
called a \emph{weakly special} subvariety of $S$ for $\VV$ (or $\Phi$).

\begin{theor}[\cite{KO}] \label{weakly special}
Let $\Phi: S^\an \to \Gamma \backslash D$ be a period map.
The bi-algebraic subvarieties of $S$ for the bi-algebraic structure
defined by $\Phi$ are precisely the weakly special subvarieties of $S$ for $\Phi$.
In analogy with \Cref{special}, they are also the closed irreducible 
algebraic subvarieties $Y \subset S$ maximal among the closed
irreducible algebraic subvarieties $Z$ of $S$ whose algebraic
monodromy group $\HH_Z$ equals $\HH_Y$.
\end{theor}

When $S=\Sh$ is a Shimura variety, these results are due to Moonen
\cite{MoMo} and \cite{UY1}. In that case the weakly special
subvarieties are also the irreducible algebraic subvarieties of $\Sh$ whose smooth locus is totally geodesic in
$\Sh^\an$ for the canonical K\"ahler-Einstein metric on $\Sh^\an=
\Gamma \backslash D$ coming
from the Bergman metric on $D$, see \cite{MoMo}.

\medskip
To study not only functional transcendence but also arithmetic
transcendence, we enrich bi-algebraic structures over $\oQ$.
A $\oQ$-bi-algebraic structure on a quasi-projective variety
$S$ defined over $\oQ$ is a bi-algebraic structure $(f: \wti{S^\an} \to Z^\an, h:
\pi_1(S^\an) \to \Aut(Z))$ such that $Z$ is
defined over $\overline{\QQ}$ and the
homomorphism $h$ takes values in $\Aut_{\oQ} Z$. An algebraic
subvariety $Y \subset \wti{S^\an}$ is
said to be defined over $\ol{\QQ}$ if its model $\ol{f(Y)}^\Zar \subset
Z$ is. A $\oQ$-bi-algebraic subvariety $W\subset S$ is an algebraic
subvariety of $S$ defined over $\oQ$ and such that any (equivalently one) of
the analytic irreducible components of $p^{-1}(W)$ is an algebraic
subvariety of $\wti{S^\an}$ defined over $\oQ$. A
$\oQ$-bi-algebraic point $s \in S(\oQ)$ is also called an {\em
  arithmetic point}. \Cref{biExample}a) is naturally defined over
$\oQ$, with arithmetic points the torsion points of $(\CC^*)^n$. In
\Cref{biExample}b) the bi-algebraic structure can be defined over $\oQ$ if
the abelian variety $A$ has CM, and its arithmetic points
are its torsion points. \Cref{biExample}c) is naturally a $\oQ$-bi-algebraic
structure, with arithmetic points the
{\em special points} of the Shimura variety (namely the special
subvarieties of dimension zero), at least when the pure part of the Shimura
variety is of Abelian type, see \cite{SW}. In all
these cases it is interesting to notice that the $\oQ$-bi-algebraic subvarieties are
the bi-algebraic subvarieties containing one arithmetic point (in
\Cref{biExample}c) these are the special subvarieties of the Shimura
variety).

The bi-algebraic structure associated with a period map $\Phi: S^\an \to \Gamma
\backslash D$ is defined over $\oQ$ as soon as $S$ is. In this case,
we expect the $\oQ$-bi-algebraic subvarieties to be precisely the special
subvarieties, see \cite[2.6 and 3.4]{K17}.

\subsection{The Ax-Schanuel theorem for period maps}
The geometry of bi-algebraic structures is controlled by the
following functional transcendence heuristic, whose idea was
introduced by Pila in the case of Shimura varieties, see\cite{Pil14},
\cite{Pil15}~:

\smallskip
\noi
{\it Ax-Schanuel principle: \label{Ax-Schanuel}
Let $S$ be an irreducible algebraic variety endowed with a non-trivial
bi-algebraic structure. Let $U \subset \wti{S^\an} \times S^\an$ be
an algebraic subvariety (for the product bi-algebraic structure) and let $W$ be an analytic irreducible component of $U
\cap \Delta$, where $\Delta$ denotes the graph of $p: \wti{S^\an} \to
S^\an$.
Then $\codim_U W \geq \dim \overline{W}^{\textnormal{bi}}$, where
$\overline{W}^{\textnormal{bi}}$ denotes the smallest bi-algebraic
subvariety of $S$ containing $p(W)$.}

\medskip
\noi
When applied to a subvariety $U \subset \wti{S^\an} \times S^\an$ of the
form $Y \times \overline{p(Y)}^\Zar$ for $Y \subset \wti{S^\an}$
algebraic, the Ax-Schanuel principle specializes
to the following:

\medskip
\noi
{\it Ax-Lindemann principle:  
Let $S$ be an irreducible algebraic variety endowed with a non-trivial
bi-algebraic structure. Let $Y\subset \wti{S^\an}$ be an
algebraic subvariety. Then $\overline{p(Y)}^\Zar$ is a bi-algebraic
subvariety of $S$.}

\medskip
Ax \cite{Ax71}, \cite{Ax72} showed that the abstract Ax-Schanuel principle holds true
for \Cref{biExample}a) and \Cref{biExample}b) above, using
differential algebra. Notice that the Ax-Lindemann
principle in \Cref{biExample}a) is the functional analog of the classical
Lindemann theorem stating that if $\alpha_1,
\dots, \alpha_n$ are $\QQ$-linearly independent algebraic numbers
then $e^{\alpha_1}, \dots, e^{\alpha_n}$ are algebraically
independent over $\QQ$. This explains the terminology.
The Ax-Lindemann principle in \Cref{biExample}c) was proven by Pila \cite{Pil} when $S$ is a product
$Y(1)^n \times (\CC^*)^k$, by Ullmo-Yafaev \cite{UY2} for projective
Shimura varieties, by Pila-Tsimerman \cite{PT} for $\cA_g$, and by
Klingler-Ullmo-Yafaev \cite{KUY} for any pure Shimura variety. The
full Ax-Schanuel principle was proven by Mok-Pila-Tsimerman for pure
Shimura varieties \cite{MPT19}.

\medskip
We conjectured in \cite[Conj. 7.5]{K17} that the Ax-Schanuel principle
holds true for the bi-algebraic structure associated to a (graded-)polarized
variation of (mixed) $\ZZ$HS on an arbitrary
quasi-projective variety $S$. Bakker and Tsimerman
proved this conjecture in the pure case:

\begin{theor}[Ax-Schanuel for $\ZZ$VHS, \cite{BT19}]\label{astheorem}
 Let $\Phi: S^\an \to \Gamma \backslash D$ be a period map.
Let $V \subset S \times \check{D}$ be an algebraic subvariety. Let $U$ be an
irreducible complex analytic component of $W \cap (S \times_{\Gamma \backslash D} D)$
such that 
\begin{equation} \label{ineg}
\codim_{S \times D} U< \codim_{S \times \check{D}} W+ \codim_{S
  \times D} (S \times_{ \Gamma \backslash D}
D)\;\;.  
\end{equation}
Then the projection of $U$ to $S$ is contained in a strict weakly
special subvariety of $S$ for $\Phi$.
\end{theor}

\begin{rem}
\cite{MPT19} was extended by Gao \cite{Gao18} to mixed
Shimura varieties of Kuga type. Recently the full Ax-Schanuel
\cite[Conj. 7.5]{K17} for variations of mixed Hodge structures has been
fully proven independently in \cite{GK} and \cite{Chiu}.
\end{rem}

The proof of \Cref{astheorem} follows a strategy
started in \cite{KUY} and fully developed in \cite{MPT19} in the
Shimura case, see \cite{Tsib} for an introduction. It does not use
\Cref{definPeriodMap} but only a week version equivalent to the
Nilpotent Orbit Theorem, and relies crucially on the definable
Chow \Cref{PS}, the Pila-Wilkie \Cref{PW} and the proof that the volume
(for the natural metric on $\Gamma \backslash D$) of the intersection
of a ball of radius $R$ in $\Gamma \backslash D$ with the horizontal
complex analytic subvariety $\Phi(S^\an)$ grows exponentially with $R$
(a negative curvature property of the horizontal tangent bundle).

\subsection{On the distribution of the Hodge locus}
\Cref{astheorem} is most useful, even in its simplest version of the
Ax-Lindemann theorem. After \Cref{CDK} one would like to understand the distribution in $S$ of the special
subvarieties for $\VV$. For instance, are there any geometric
constraints on the Zariski closure of $\HL(S, \VV^\otimes)$? To
approach this question, let us decompose the adjoint group $\GG^\ad$
into a product $\GG_1 \times \cdots \times \GG_r$ of its simple factors. It gives rise
(after passing to a finite \'etale covering if necessary) to a
decomposition of the Hodge variety $\Gamma \backslash D$ into a
product of Hodge varieties $\Gamma_1\backslash D_1 \times \cdots \times \Gamma_r
\backslash D_r$. A special subvariety $Z$ of $S$ for $\VV$ is said of \emph{positive
period dimension} if $\dim_\CC \Phi(Z^\an) >0$; and of
\emph{factorwise positive period dimension} if moreover the projection
of $\Phi(Z^\an)$ on each factor $\Gamma_i \backslash D_i$ has positive
dimension. The {\em Hodge locus of factorwise positive period dimension} 
  $\HL(S, \VV^\otimes)_{\fpos}$ is the union of the strict special
subvarieties of positive period dimension, it is contained
in the \emph{Hodge locus of positive period dimension} $\HL(S, \VV^\otimes)_\pos$ union of the strict special
subvarieties of positive period dimension, and the two coincide if
$\GG^\ad$ is simple.

\medskip
Using the Ax-Lindemann theorem special case of \Cref{astheorem} and a global
algebraicity result in the total bundle of $\cV$, Otwinowska and the
author proved the following:

\begin{theor}[\cite{KO}] \label{KO}
Let $\VV$ be a polarized $\ZZ$VHS on a smooth connected complex quasi-projective variety
$S$. Then either $\HL(S, \VV^{\otimes})_\fpos$ is Zariski-dense in
$S$; or it is an algebraic subvariety of $S$ (i.e., the set of strict special subvarieties of $S$ for 
$\VV$ of factorwise positive period dimension has only finitely many maximal
elements for the inclusion).
\end{theor}

\begin{Example} \label{shimura intersection}
The simplest example of \Cref{KO} is the following. Let $S \subset \cA_g$ be a Hodge-generic closed irreducible
subvariety. Either the set of positive dimensional closed irreducible subvarieties
of $S$ which are not Hodge generic has 
finitely many maximal elements (for the inclusion), or their union is Zariski-dense in
$S$.
\end{Example}
 
\begin{Example} \label{Noether-Lefschetz}
Let $B \subset \proj H^0(\proj^3_\CC, \cO(d))$ be the open subvariety
parametrising the smooth surfaces of degree $d$ in $\proj^3_\CC$. Suppose
$d>3$. The classical Noether theorem
states that any surface $Y \subset \proj^3_\CC$ corresponding to a
very general point $[Y] \in B$ has Picard group $\ZZ$:
every curve on $Y$ is a complete intersection of $Y$ with another
surface in $\proj^3_\CC$. The countable union $\NL(B)$ of closed algebraic
subvarieties of $B$ corresponding to surfaces with bigger Picard group
is called the Noether-Lefschetz locus of $B$.
Let $\VV \to B$ be the $\ZZ$VHS $R^2f_*\ZZ_\prim$, where $f: \cY \to
B$ denotes the universal family of surfaces of degree $d$. Clearly $\NL(B) \subset \HL(B,
\VV^\otimes)$. Green (see \cite[Prop.5.20]{Voisin02}) proved
that $\NL(B)$, hence also $\HL(B, \VV^\otimes)$, is analytically dense in $B$. 
Now \Cref{KO} implies the following:
Let $S \subset B$ be a Hodge-generic closed irreducible subvariety.
Either $S \cap \HL(B, \VV^\otimes)_\fpos$ contains only finitely many maximal positive
dimensional closed irreducible subvarieties of $S$, or the union of such subvarieties is Zariski-dense in $S$.
\end{Example}

\section{Typical and atypical intersections: the Zilber-Pink conjecture for period
maps} \label{atypicalsection}

\subsection{The Zilber-Pink conjecture for $\ZZ$VHS: conjectures}

In the same way that the Ax-Schanuel principle controls the geometry of
bi-algebraic structures, the diophantine geometry of
$\oQ$-bi-algebraic structures is controlled by the following heuristic:

\medskip
\noi
{\it Atypical intersection principle: Let $S$ be an irreducible
  algebraic $\oQ$-variety endowed with a 
  $\oQ$-bi-algebraic structure. Then the union $S_\atyp$ of atypical
  $\oQ$-bi-algebraic subvarieties of $S$ is an algebraic subvariety of $S$
  (i.e. it contains only finitely many atypical $\oQ$-bi-algebraic
  subvarieties maximal for the inclusion).}

\medskip
\noi
Here a $\oQ$-bi-algebraic subvariety $Y \subset S$ is said to be {\em atypical} for
the given bi-algebraic structure on $S$ if it is obtained as an excess
intersection of $f(\wti{S^\an})$ with its model
$\ol{f(\tilde{Y})}^\Zar \subset Z$; and $S_\atyp$ denotes
the union of all atypical subvarieties of $S$.
As a particular case of the atypical intersection principle:

\medskip
\noi
{\it Sparsity of arithmetic points principle: Let $S$ be an irreducible
  algebraic $\oQ$-variety endowed with a 
  $\oQ$-bi-algebraic structure. Then any irreducible algebraic subvariety of $S$
  containing a Zariski-dense set of atypical arithmetic points is a
  $\oQ$-bi-algebraic subvariety.}

\medskip
This principle that arithmetic points are sparse is a theorem of Mann
\cite{Ma65} in \Cref{biExample}a). For abelian varieties over $\oQ$ (\Cref{biExample}b)) this is
the Manin-Mumford conjecture proven first by Raynaud \cite{Ray},
saying that an irreducible subvariety of an abelian variety over $\oQ$
containing a Zariski-dense set of torsion point is the translate of an
abelian subvariety by a torsion point. For Shimura varieties of
abelian type (\Cref{biExample}c)) this is the classical Andr\'e-Oort conjecture
\cite{Andre}, \cite{Oort} stating that an irreducible subvariety of a
Shimura variety containing a Zariski-dense set of special points is special. It
has been proven in this case using tame geometry and following the strategy proposed by
Pila-Zannier \cite{PZ} (let us mention \cite{Pil}, \cite{U2},
\cite{PT}, \cite{KUY}, \cite{AGHM}, \cite{YuZh}, \cite{Tsi}; and
\cite{Gao} in the mixed case; see
\cite{KUY2} for a survey). Recently the Andr\'e-Oort conjecture in
full generality has been obtained in
\cite{PST}, reducing to the case of abelian type using ingredients
from $p$-adic Hodge theory. We refer to \cite{Zannier} for many
examples of atypical intersection problems.

\medskip
In the case of Shimura varieties (\Cref{biExample}c)) the general atypical intersection
principle is the Zilber-Pink conjecture \cite{Pink05}, \cite{Zil02},
\cite{panor}. Only very few instances of the Zilber-Pink conjecture
are known outside of the Andr\'e-Oort conjecture, see \cite{HP12},
\cite{HP16}, \cite{DR} for example.

\medskip
For a general polarized $\ZZ$VHS $\VV$ with period map $\Phi: S^\an
\to \Gamma \backslash D$, which we can assume to be proper without
loss of generality, we already mentioned that even the geometric characterisation of the $\oQ$-bi-algebraic
subvarieties as the special subvarieties is unknown. Replacing the $\oQ$-bi-algebraic subvarieties of $S$ by the
special ones, we define:

\begin{defi} \label{atypical}
A special subvariety $Z =
\Phi^{-1}(\Gamma_Z \backslash D_Z)^0\subset S$ is said \emph{atypical} if
either $Z$ is {\em singular for $\VV$} (meaning that $\Phi(Z^\an)$ is
contained in the singular locus of the complex analytic variety $\Phi(S^\an)$), or if $\Phi(S^{\an})$ and $\Gamma_Z
  \backslash D_Z$ do not intersect generically along $\Phi(Z)$: 
  $$ \codim_{\Gamma\backslash D} \Phi(Z^{\an}) < \codim_{\Gamma\backslash
    D} \Phi(S^{\an}) + \codim_{\Gamma\backslash D} \Gamma_Z \backslash
  D_Z\;\;.$$
  \noindent
  Otherwise it is said to be \emph{typical}.
  \end{defi}

Defining the \emph{atypical Hodge locus} $\HL(S,\VV^\otimes)_{\atyp}
\subset \HL(S, \VV^\otimes)$ as the union of the atypical special
subvarieties of $S$ for $\VV$, we obtain the following precise atypical
intersection principle for $\ZZ$VHS, first proposed in \cite{K17} in a more
restrictive form:

\begin{conj}[Zilber--Pink conjecture for $\ZZ$VHS, \cite{K17}, \cite{BKU}] \label{main conj}
  Let $\VV$ be a polarizable $\ZZ$VHS on an irreducible smooth
  quasi-projective variety $S$. The atypical Hodge locus
  $\HL(S,\VV^\otimes)_{\atyp}$ is a finite union of atypical special
subvarieties of $S$ for $\VV$. Equivalently: the set of atypical
special subvarieties of $S$ for $\VV$ has finitely many maximal elements for
the inclusion.
\end{conj}

Notice that this conjecture is in some sense more general than the above
atypical intersection principle, as we don't assume that $S$ is defined over $\oQ$; this has to
be compared to the fact that the Manin-Mumford conjecture holds true
for every complex abelian variety, not necessarily defined over
$\oQ$.

\begin{Example} \label{BaldiUllmo}
Recently Baldi and Ullmo \cite{BU} proved a special case of \Cref{main conj} of
much interest. Margulis' arithmeticity theorem states that any lattice
in a simple real Lie group $G$ of real rank at least $2$ is arithmetic: it is commensurable with a group
$\GG(\ZZ)$, for $\GG$ a $\QQ$-algebraic group such that $\GG(\RR)= G$
up to a compact factor. On the other hand the structure of lattices in
a simple real Lie group of rank $1$, like the group $PU(n,1)$ of
holomorphic isometries of the complex unit ball $\BB^n_\CC$ endowed
with its Bergman metric, is an essentially open
question. In particular there exist non-arithmetic lattices in
$PU(n,1)$, $n=2,3$. Let $\iota: \Lambda \hookrightarrow PU(n,1)$ be a lattice. The ball
quotient $S^\an:= \Lambda \backslash \mathbf{B}^n_\CC$ is the
analytification of a complex algebraic variety $S$. By results of
Simpson and Esnault-Groechenig, there exists a $\ZZ$VHS $\Phi: S^\an
\to \Gamma \backslash (\BB^n_\CC \times D')$ with monodromy
representation $\rho : \Lambda \to PU(n,1) \times G'$ whose first
factor $\Lambda \to PU(n,1)$ is the rigid representation
$\iota$. The special subvarieties of $S$ for $\VV$ are the totally
geodesic complex subvarieties of $S^\an$. When $\Lambda$ is
non-arithmetic, they are automatically atypical. In accordance with
\Cref{main conj} in this case, Baldi and Ullmo prove that if $\Lambda$ is non-arithmetic,
then $S^\an$ contains only finitely many maximal totally geodesic 
subvarieties. This result has been proved independently by Bader,
Fisher, Miller, and Stover \cite{BFMS}, using completely different
methods from homogeneous dynamics.
\end{Example}

Among the special points for a $\ZZ$VHS $\VV$, the CM-points
(i.e. those for which the Mumford-Tate
group is a torus) are always atypical except if the generic Hodge
datum $(\GG, D)$ is of Shimura type and the period map $\Phi$ is
dominant. Hence, as explained in \cite[Section
5.2]{K17}, \Cref{main conj} implies the following:

\begin{conj}[Andr\'e-Oort conjecture for $\ZZ$VHS, \cite{K17}] \label{AO}
 Let $\VV$ be a polarizable $\ZZ$VHS on an irreducible smooth
  quasi-projective variety $S$. If $S$ contains a Zariski-dense set of
  CM-points then the generic Hodge datum $(\GG, D)$ of $\VV$ is a
  Shimura datum, and the period map $\Phi: S^\an \to \Gamma \backslash
  D$ is an algebraic map, dominant on the Shimura variety $\Gamma \backslash
  D$.
  \end{conj}

  \begin{Example}
Consider the Calabi-Yau Hodge structure $V$ of weight $3$ with Hodge numbers
$h^{3,0}= h^{2,1}=1$ given by the mirror dual
quintic. Its universal deformation space $S$ is the projective line minus
3 points, which carries a $\ZZ$VHS $\VV$ of the same type. This gives a
non-trivial period map $ \Phi: S^{\an} \to \Gamma \backslash D$, where $D=
\mathbf{Sp}(4, \RR)/U(1)\times U(1)$ is a $4$-dimensional period domain. This period map is known not to
factorize through a Shimura subvariety (its algebraic monodromy group
is $\Sp_4$). \Cref{AO} in that case predicts that $S$
contains only finitely many points CM-points $s$. A version of this prediction
already appears in \cite{GuVa}. The more general \Cref{main conj} also
predicts that $S$ contains only finitely many points $s$
where $\VV_s$ splits as a direct sum of two (Tate twisted) weight one
Hodge structures $(\VV_s^{2,1}
\oplus \VV_s^{1,2})$ and its orthogonal for the Hodge metric $ (\VV_s^{3,0}
\oplus \VV_s^{0,3})$ (the so-called ``rank two attractors'' points, see
\cite{Moore}).
\end{Example}
 
\Cref{main conj} about the atypical Hodge locus takes all its meaning if we compare it to the
expected behavior of its complement, the {\em typical Hodge locus}
$\HL(S, \VV^\otimes)_\typ:= \HL(S, \VV^\otimes) \setminus \HL(S, \VV^\otimes_\atyp)$:

\begin{conj}[Density of the typical Hodge locus, \cite{BKU}] \label{conj-typical}
If $\HL(S, \VV^\otimes)_\typ$ is not empty then it is dense (for the analytic topology) in
$S^\an$.
\end{conj}

\Cref{main conj} and \Cref{conj-typical} imply immediately the
following, which clarifies the possible alternatives in \Cref{KO}:

\begin{conj}[\cite{BKU}] \label{second-main}
  Let $\VV$ be a polarizable $\ZZ$VHS on an irreducible smooth
  quasi-projective variety $S$. If $\HL(S, \VV^\otimes)_\typ$ is
  empty then $\HL(S, \VV^\otimes)$ is algebraic; otherwise
  $\HL(S, \VV^\otimes)$ is analytically dense in $S^\an$.
  \end{conj}

\subsection{The Zilber-Pink conjecture for $\ZZ$VHS: results}
In \cite{BKU} Baldi, Ullmo and I establish the \emph{geometric part} of
\Cref{main conj}: the maximal atypical special subvarieties {\em of positive
period dimension} arise in a finite number of families whose geometry
is well-understood. We can't say anything on
the atypical locus of zero period dimension (for which different ideas
are certainly needed):

\begin{theor}[Geometric Zilber--Pink, \cite{BKU}]\label{geometricZP}
Let $\VV$ be a polarizable $\ZZ$VHS on a smooth connected complex quasi-projective variety
$S$. Let $Z$ be an irreducible component of the Zariski closure
of $\HL(S, \VV^{\otimes})_{\pos, \atyp}:= \HL(S,
\VV^\otimes)_{\pos} \cap \HL(S,
\VV^\otimes)_{\atyp}$ in $S$. Then:

\begin{itemize}
  \item[(a)] Either $Z$
    is a maximal atypical special subvariety;
    \item[(b)] 
Or the generic adjoint Hodge datum $(\GG_Z^\ad, D_{G_{Z}})$
  decomposes as a non-trivial product $(\GG', D') \times (\GG'', D'')$,
  inducing (after replacing $S$ by a finite \'{e}tale cover if
  necessary)
  \begin{equation*}
\Phi_{|Z^\an}= (\Phi', \Phi''): Z^\an \to  \Gamma_{\GG_{Z}}\backslash D_{G_{Z}}= \Gamma'
\backslash D' \times  \Gamma'' \backslash D''\subset \Gamma \backslash
D,
\end{equation*}
such that $Z$ contains a Zariski-dense
set of atypical special subvarieties for $\Phi''$ of zero
period dimension. Moreover $Z$ is Hodge generic in the special subvariety $\Phi^{-1}(
\Gamma_{\GG_{Z}}\backslash D_{G_{Z}})^0$ of $S$ for $\Phi$, which is typical. 
 \end{itemize}
\end{theor}

\noindent
\Cref{main conj}, which also takes into account the atypical special
subvarieties of zero period dimension, predicts that the branch (b) of the
  alternative in the conclusion of \Cref{geometricZP} never
  occurs. \Cref{geometricZP} is proven using properties of definable
  sets and the Ax-Schanuel \Cref{astheorem}, following an idea
  originating in \cite{U2}.

  \medskip
As an application of \Cref{geometricZP}, let us consider the
\emph{Shimura locus} of $S$ for $\VV$, namely the
union of the special subvarieties of $S$ for $\VV$ which are of Shimura
type (but not necessarily with dominant period maps). In
\cite{K17}, I asked (generalizing the Andr\'e-Oort conjecture
for $\ZZ$VHS) whether a polarizable $\ZZ$VHS
$\VV$ on $S$ whose Shimura locus in Zariski-dense in $S$ is
necessarily of Shimura type. As a corollary of \Cref{geometricZP} we obtain:

\begin{theor}[\cite{BKU}] \label{shimuralocus}
Let $\VV$ be a polarizable $\ZZ$VHS on a smooth irreducible complex quasi-projective variety
$S$, with generic Hodge datum $(\GG, D)$. Suppose that the Shimura
locus of $S$ for $\VV$ \emph{of positive period dimension} is Zariski-dense in $S$.
If $\GG^\ad$ is simple then $\VV$ is of Shimura
type.
\end{theor}

\subsection{On the algebraicity of the Hodge
  locus}

In view of \Cref{second-main}, it is natural to ask if there a simple combinatorial criterion on $(\GG, D)$ for deciding
  whether $\HL(S, \VV)_\typ$ is empty. Intuitively, one expects that
the more ``complicated'' the Hodge structure is, the smaller the typical Hodge locus
should be, due to the constraint imposed by Griffiths'
transversality. Let us measure the complexity of $\VV$ by its {\em
  level}: when $\GG^\ad$ is simple, it is the greatest integer $k$ such that $\Fg^{k,
  -k}\not =0$ in the Hodge decomposition of the Lie algebra $\Fg$ of
$\GG$; in general one takes the minimum of these integers obtained for
each simple $\QQ$-factor of $\GG^\ad$. While strict typical special subvarieties usually abound for   
$\ZZ$VHSs of level one (e.g. families of abelian varieties, see \Cref{shimura intersection}; or families
of K3 surfaces) and can occur in level two (see
\Cref{Noether-Lefschetz}), they do not exist in level at least three!

\begin{theor}[\cite{BKU}]\label{level criterion}
Let $\VV$ be a polarizable $\ZZ$VHS
on a smooth connected complex quasi-projective variety 
$S$. If $\VV$ is of level at
 least $3$ then $\HL(S,\VV^\otimes)_{\typ} = \emptyset$ (and thus $\HL(S,
 \VV^\otimes)= \HL(S, \VV^\otimes)_\atyp$).
\end{theor}

The proof of \Cref{level criterion} is purely Lie theoretic. Let
$(\GG, D)$ be the generic Hodge datum of $\VV$ and $\Phi: S^\an \to
\Gamma \backslash D$ its period map. Suppose
that $Y \subset S$ is a typical special subvariety, with
generic Hodge datum $(\GG_Y, D_Y)$. The typicality condition and the
horizontality of the period map $\Phi$ imply that $\Fg_Y^{-i, i}=
\Fg^{-i,i}$ for all $i \geq 2$ (for the Hodge structures on the Lie
algebras $\Fg_Y$ and $\Fg$ defined by some point of $D_Y$). Under the assumption that $\VV$
has level at least $3$, we show that this is enough to ensure that $\Fg_Y =
\Fg$, hence $Y=S$. Hence there are no strict typical special subvariety.

\medskip
Notice that \Cref{main conj} and \Cref{level criterion} imply:

\begin{conj}[Algebraicity of the Hodge locus in level at least 3,
  \cite{BKU}] \label{conjAlg}
Let $\VV$ be a polarizable $\ZZ$VHS
on a smooth connected complex quasi-projective variety 
$S$. If $\VV$ is of level at least $3$ then $\HL(S, \VV^\otimes)$ is
algebraic.
\end{conj}

The main result of \cite{BKU}, which follows immediately from \Cref{geometricZP} and
\Cref{level criterion}, is the following stunning geometric
reinforcement of \Cref{CDK} and \Cref{KO}:

\begin{theor}[\cite{BKU}] \label{corol}
If $\VV$ is of level at least $3$ then $\HL(S,
\VV^\otimes)_\fpos$ is algebraic. 
\end{theor}

As a simple geometric illustration of \Cref{corol}, we prove the
following, to be contrasted with the $n=2$ case (see \Cref{Noether-Lefschetz}):
\begin{cor} \label{hypersurface}
Let $\PP^{N(n, d)}_\CC$ be the projective space parametrising the hypersurfaces $X$ of
$\PP^{n+1}_\CC$ of degree $d$ (where $N(n, d)=\binom{n+d+1}{d}-1$). Let $U_{n, d} \subset \PP^{N(n, d)}_\CC$ be
the Zariski-open subset parametrising the smooth
hypersurfaces $X$ and let 
$
\VV \to U_{n,d}
$
be the $\ZZ$VHS corresponding to the primitive cohomology
$H^n(X, \ZZ)_\prim$. If $n\geq 3$ and $d>5$ then $\HL(U_{n,d}, \VV^\otimes)_\pos \subset U_{n,d}$ is algebraic.
\end{cor}

\subsection{On the typical Hodge locus in level one and two}

In the direction of \Cref{conj-typical}, we obtain:

\begin{theor}[Density of the typical locus, \cite{BKU}] \label{typicallocus}
Let $\VV$ be a polarized $\ZZ$VHS
on a smooth connected complex quasi-projective variety 
$S$. If the typical Hodge locus $\HL(S,\VV^\otimes)_{\typ}$ is
non-empty (hence the level of $VV$ is one or two by \Cref{level
  criterion}) then
$\HL(S,\VV^\otimes)$ is analytically (hence Zariski) dense in
$S$.
\end{theor}

\noi
Notice that, in \Cref{typicallocus}, we also treat
the typical Hodge locus of zero period dimension. 
\Cref{typicallocus} is new even for $S$ a subvariety of a Shimura variety. Its proof is inspired by the arguments of Chai
\cite{chai} in that case.

\medskip
It remains to find a criterion for deciding whether, in level one or
two, the typical Hodge locus $\HL(S,\VV^\otimes)_{\typ} $ is empty or
not. We refer to \cite[Theor. 2.15]{KOU} and \cite{Tayou},
\cite{TayouTholozan} for results in this direction.

\section{Arithmetic aspects} \label{arithm}

We turn briefly to some arithmetic aspects of period maps.

\subsection{Field of definition of special subvarieties} \label{fields}
Once more the geometric case provides us with a motivation and a
heuristic. Let $f:X \to S$ be a smooth projective morphism of
connected algebraic varieties defined 
over a number field $L\subset \CC$ and let $\VV$ be the natural polarizable
$\ZZ$VHS on $S^\an$ with underlying local system $R^\bullet f_*^\an \ZZ$. In that case, the Hodge conjecture
implies that each special subvariety $Y$ of $S$ for $\VV$ is defined
over $\ol{\QQ}$ and that each of the $\Gal(\oQ/L)$-conjugates of $Y$
is again a special subvariety of $S$ for $\VV$. More generally, let us
say that a polarized
$\ZZ$VHS $\VV= (\VV_\ZZ, (\cV, F^\bullet,\nabla), 
\mathrm{q})$ on $S^\an$ is {\em defined over a number field $L \subset
  \CC$} if $S$, $\cV$, $F^\bullet$ and $\nabla$ are defined over 
$L$ (with the obvious compatibilities).

\begin{conj} \label{conj1}
  Let $\VV$ be a $\ZZ$VHS defined over a number field $L \subset
  \CC$. Then any special subvariety of $S$ for $\VV$ is
defined over $\oQ$, and any of its finitely many $\Gal(\oQ/L)$-conjugates is a
  special subvariety of $S$ for $\VV$.
\end{conj}

There are only few results in that
direction: see \cite[Theor. 0.6]{V07} for a proof under a strong
geometric assumption; and \cite{SS}, where it is shown that
when $S$ (not necessarily $\VV$) is defined over $\oQ$, then a special subvariety of $S$ for
$\VV$ is defined over $\oQ$ if and only if it contains a $\oQ$-point
of $S$. In \cite{KOU} Otwinowska, Urbanik and I provide a simple geometric
criterion for a special subvariety of $S$ for $\VV$ to satisfy
\Cref{conj1}. In particular we obtain:

\begin{theor}[\cite{KOU}] \label{mainKOU}
Let $\VV$ be a polarized $\ZZ$VHS on a smooth connected complex
quasi-projective variety $S$. Suppose that the adjoint generic Mumford-Tate group
$\GG^\ad$ of $\VV$ is simple. If $S$ is defined over a number field $L$,  
then any maximal (strict) special subvariety $Y \subset S$ of positive
period dimension is defined over $\oQ$. If moreover $\VV$ is defined over $L$ then the finitely many
      $\Gal(\oQ/L)$-translates of $Y$ are special
      subvarieties of $S$ for $\VV$.
    \end{theor}

As a corollary of \Cref{corol} and \Cref{mainKOU}, one obtains the
following, which applies for instance in the situation of
\Cref{hypersurface}.

\begin{cor} \label{corolKOU}
Let $\VV$ be a polarized variation of $\ZZ$-Hodge structure on a
smooth connected 
quasi-projective variety $S$. Suppose that $\VV$ is of level at
least $3$, and that it is defined over $\oQ$. Then $\HL(S,
\VV^\otimes)_\fpos$ is an algebraic subvariety of $S$, defined over
$\oQ$.
\end{cor}

It is interesting to notice that \Cref{conjAlg}, which is stronger
than \Cref{corol}, predicts the existence of a Hodge generic
$\oQ$-point in $S$ for $\VV$ in the situation of \Cref{corolKOU}.

\medskip
As the criterion given in \cite{KOU} is purely geometric, it says
nothing about fields of definitions of special points. It is however
strong enough to reduce the first part of \Cref{conj1} to this particular
case:
\begin{theor} \label{thm2} ~
Special subvarieties for $\ZZ$VHSs defined over $\oQ$ are
defined over $\oQ$ if and only if it holds true for special
points.
\end{theor}

\subsection{Absolute Hodge locus} \label{absoluteSection}
Interestingly, \Cref{conj1} in the geometric case follows from an
\emph{a priori} much
weaker conjecture than the Hodge conjecture. Let $f: X\to S$ be a
smooth projective morphism of smooth connected complex quasi-projective
varieties. For any automorphism $\sigma \in \Aut\, (\CC/\QQ)$ we can
consider the algebraic family $f^\sigma: X^\sigma \to S^\sigma$, where
$\sigma^{-1}: S^\sigma= S\times_{\CC, \sigma} \CC
\xrightarrow{\sim} S$ is the natural
isomorphism of abstract schemes; and the attached polarizable
$\ZZ$VSH $\VV^\sigma= (\VV^\sigma_\ZZ, \cV^\sigma, {F^\bullet}^\sigma,
\nabla^\sigma)$ with underlying local system $\VV^\sigma_\ZZ =
R{f^\sigma}_* ^\an \ZZ$ on $(S^\sigma)^\an$. The algebraic
construction of the algebraic de Rham cohomology provides compatible canonical
comparison isomorphisms $\iota^\sigma: (\cV^\sigma,
{F^\bullet}^\sigma, \nabla^\sigma) \xrightarrow{\sim} {\sigma^{-1}}^* (\cV,
F^\bullet, \nabla)$ of the associated algebraic filtered vector
bundles with connection. More generally a collection of $\ZZ$VHS
$(\VV^\sigma)_\sigma$ with such compatible comparison isomorphisms is called a {\em (de Rham)
  motivic variation of Hodge structures} on $S$, in which case we
write $\VV:= \VV^{\Id}$. Following Deligne (see \cite{CS14} for a nice
exposition), an {\em
  absolute Hodge tensor} for such a collection is a Hodge tensor
  $\alpha$ for $\VV_s$ such that the conjugates ${\sigma^{-1}}^*
  \alpha_{\dR}$ of the De Rham component of $\alpha$ defines a Hodge
  tensor in $\VV^\sigma_{\sigma(s)}$ for all $\sigma$. The
  \emph{generic absolute Mumford-Tate group} for $(\VV^\sigma)_\sigma$
  is defined in terms of the absolute Hodge tensors as the generic
  Mumford-Tate group is defined in terms of the Hodge tensors. Thus $\GG
  \subset \GG^\AH$. In view of \Cref{special} the following is
  natural:

  \begin{defi} \label{absolute special}
    Let $(\VV^\sigma)_\sigma$ be a (de Rham) motivic variation of
    Hodge structure on a smooth connected complex quasi-projective
    variety $S$. A closed irreducible algebraic subvariety
    $Y$ of $S$ is called {\em absolutely special} if it is maximal
    among the closed irreducible algebraic subvarieties $Z$ of $S$
    satisfying $\GG_Z^\AH= \GG_Y^\AH$.
  \end{defi}
  
  In the geometric case, the Hodge conjecture implies,
  since any automorphism $\sigma \in \Aut\, (\CC/\QQ)$ maps algebraic cycles in
  $X$ to algebraic cycles on $X^\sigma$, 
  the following conjecture of Deligne:

  \begin{conj}[\cite{Del82}] \label{absolute}
    Let $(\VV^\sigma)_\sigma$ be a (de Rham) motivic variation of
    Hodge structure on $S$. Then all Hodge tensors are absolute Hodge
    tensors, i.e. $\GG= \GG^{\AH}$.
  \end{conj}

  This conjecture immediately implies:

  \begin{conj} \label{conjabsspecial}
   Let $(\VV^\sigma)_\sigma$ be a (de Rham) motivic variation of
    Hodge structure on $S$. Then any special subvariety of $S$ for
    $\VV$ is absolutely special for $(\VV^\sigma)_\sigma$.
\end{conj}

Let us say that a (de Rham) motivic variation $(\VV^\sigma)_\sigma$ is defined over
$\oQ$ if $\VV^\sigma = \VV$ for all $\sigma \in \Aut
(\CC/\ol{\QQ})$. In the geometric case any morphism $f: X
\to S$ defined over $\oQ$ defines such a (de Rham) motivic variation $(\VV^\sigma)_\sigma$ over
$\oQ$. Notice that the absolutely special subvarieties of $S$ for
$(\VV^\sigma)_\sigma$  are then by their very definition defined over
$\oQ$, and their Galois conjugates are also special. In particular
\Cref{conjabsspecial} implies \Cref{conj1} in the geometric
case. As proven in \cite{V07}, Deligne's conjecture is actually
equivalent to a much stronger version of \Cref{conj1}, where one
replaces the special subvarieties of $S$ (components of the Hodge locus) with the special subvarieties in
the total bundle of $\cV^\otimes$ (components of the locus of Hodge
tensors).

\medskip
Recently T.Kreutz, using the same geometric argument as in
\cite{KOU}, justified \Cref{mainKOU} by proving:

\begin{theor}[\cite{Kreutz1}]
  Let $(\VV^\sigma)_\sigma$ be a (de Rham) motivic variation of
 Hodge structure on $S$. Suppose that the adjoint generic Mumford-Tate
 group $\GG^\ad$ is simple. Then any strict maximal special subvariety
 $Y\subset S$ of positive period dimension for $\VV$ is absolutely
 special.
\end{theor}

\noindent
We refer the reader to \cite{Kreutz2}, as well as \cite{Urbanik}, for
other arithmetic aspects of Hodge loci taking into account not only
the de Rham incarnation of absolute Hodge classes but also their $\ell$-adic components.


\bigskip
\noindent Bruno Klingler : Humboldt Universit\"at zu Berlin

\noindent email : \texttt{bruno.klingler@hu-berlin.de}


\begin{thebibliography}{99}







\bibitem[An89]{Andre} Y. Andr\'e, \emph{G-functions and geometry},
  Aspects of Mathematics E13 (1989)

 \bibitem[An92]{An92} Y. Andr\'e, \emph{Mumford-Tate groups of mixed
     Hodge   
structures and the theorem of the fixed part},
Compositio Math. {\textbf 82} (1992) 1-24

\bibitem[AGHM]{AGHM} F. Andreatta, E. Goren, B. Howard,
 K. Madapusi-Pera, \emph{Faltings heights of abelian varieties with
    complex multiplication}, Ann. of Math (2) \textbf{187} (2018), no.2,
  391-531 


\bibitem[Art70]{Artin} M. Artin, \emph{Algebraization of formal
    moduli. II. Existence of modifications}, Annals of
  Math. \textbf{91}, (1970), 88--135

\bibitem[Ax71]{Ax71} J. Ax, \emph{On Schanuel's conjecture}, Annals of
  Math. \textbf{93} (1971), 1-24

\bibitem[Ax72]{Ax72} J. Ax, \emph{Some topics in differential algebraic
    geometry. I. Analytic subgroups of algebraic groups},
  Amer. J. Math. \textbf{94} (1972), 1195-1204

\bibitem[BFMS20]{BFMS} U. Bader, D. Fisher, N. Miller, and M. Stover,
  \emph{Arithmeticity, superrigidity and totally geodesic submanifolds
    of complex hyperbolic manifolds}
  \url{http://arxiv.org/abs/2006.03008}
  
\bibitem[BB66]{BB66} W.L. Baily, A. Borel, \emph{Compactification of arithmetic  
quotients of bounded symmetric domains},
Ann. of Math., \textbf{84} (1966), 442-528

\bibitem[BT19]{BT19} B. Bakker, J. Tsimerman, \emph{The Ax-Schanuel
      conjecture for variations of Hodge structures},
    Invent. Math. \textbf{217} (2019), no.1, 77-94

\bibitem[BBKT20]{BBKT} B. Bakker, Y. Brunebarbe, B. Klingler,
    J. Tsimerman, \emph{Definability of mixed period maps},
    \url{http://arxiv.org/abs/2006.12403}, to appear in  J. Eur. Math. Soc. 
    
\bibitem[BKT20]{BKT} B. Bakker, B. Klingler, J. Tsimerman, \emph{Tame
    topology of arithmetic quotients and algebraicity of Hodge loci},
  J. Amer. Math. Soc. 33 (2020), 917-939

  \bibitem[BBT18]{BBT} B. Bakker, Y. Brunebarbe, J.Tsimerman, \emph{
      o-minimal GAGA and a conjecture of Griffiths},
    \url{http://arxiv.org/abs/1811.12230}

\bibitem[BKU21]{BKU} G. Baldi, B. Klingler, E. Ullmo, \emph{On the
        distribution of the Hodge locus},
      \url{http://arxiv.org/abs/2107.08838}

\bibitem[BU20]{BU} G. Baldi, E. Ullmo, \emph{Special
subvarieties of non-arithmetic ball quotients and Hodge
theory}, \url{http://arxiv.org/abs/2005.03524}

\bibitem[BoP89]{BP} E. Bombieri, J.Pila, \emph{The number of integral
    points on arcs and ovals}, Duke Math. J. \textbf{59} (1989), no.2,
  337--357. 



  \bibitem[Bor69]{bor} A. Borel, {\it Introduction aux groupes arithm{\'e}tiques},
   Publications de l'Institut de Math{\'e}matique de l'Universit{\'e} de
   Strasbourg, XV. Actualit{\'e}s Scientifiques et Industrielles, No. 1341
   Hermann, Paris (1969)
   
  \bibitem[Bor72]{Bor72} A. Borel, {\it Some metric properties of
    arithmetic quotients of symmetric spaces and an extension
    theorem}, J. Differential Geometry {\bf 6} (1972), 543-560

  \bibitem[BP09a]{BP1} P. Brosnan, G. Pearlstein, \emph{Zero loci of
      admissible normal functions with torsion singularities}, Duke
    Math. J. \textbf{150} (2009), no. 1, 77-100
  
  \bibitem[BP09b]{BP2} P. Brosnan, G. Pearlstein, \emph{The zero locus of
    an admissible normal function}, Ann. of Math. (2) \textbf{170} (2009),
  no. 2, 883-897
  
  \bibitem[BPS10]{BPS} P. Brosnan, G. Pearlstein, C. Schnell, \emph{The
    locus of Hodge classes in an admissible variation of mixed Hodge
    structure}, C.R. Math. Acad. Sci. Paris \textbf{348} (2010),
  no. 11-12, 657-660
  
  \bibitem[BP13]{BP3} P. Brosnan, G. Pearlstein, \emph{On the
      algebraicity of the zero locus of an admissible normal
      function}, Compos. Math. \textbf{149} (2013, no. 11, 1913--1962

\bibitem[CDK95]{CDK95} E. Cattani, P. Deligne, A. Kaplan, \emph{On the locus of Hodge classes.}
J. of AMS, 8 (1995), 483-506

\bibitem[CKS86]{CKS} E. Cattani, A. Kaplan, W. Schmid, \emph{Degeneration of
    Hodge structures}, Ann. of Math. (2) \textbf{123} (1986), no. 3,
  457-535

  \bibitem[Chai98]{chai} C.L. Chai, {\em Density of members with extra
    Hodge cycles in a family of Hodge structures},  Asian
  J. Math. {\bf 2} (1998), no. 3, 405-418

  \bibitem[CS14]{CS14} F. Charles, C. Schnell, \emph{Notes on absolute
      Hodge classes}, Hodge theory, Math. Notes \textbf{49} (2014),
    469-530

\bibitem[Chiu21]{Chiu} K.C.T. Chiu, \emph{Ax-Schanuel for variations
    of mixed Hodge structures}, \url{http://arxiv.org/abs/2101.10968}

\bibitem[DR18]{DR} C. Daw, J.Ren, \emph{Applications of the hyperbolic
    Ax-Schanuel conjecture}, Compos. Math. \textbf{154} (2018) no.9,
  1843-1888
  
\bibitem[Del70]{Del} P. Deligne, \emph{Equations diff\'erentielles \`a
    points singuliers r\'eguliers}, LNM {\bf 163} (1970)
  
\bibitem[De71a]{Del71}  P. Deligne, \emph{Th\'eorie de Hodge II}, 
    Publ. Math. IHES \textbf{40} (1971) 5-57

\bibitem[De71b]{De1} P. Deligne, \emph{Travaux de Shimura}, S\'eminaire
  Bourbaki Expos\'e 389, LNM {\bf 244} (1971), 123-165


\bibitem[De72]{Del72}  P. Deligne, {\it La conjecture de Weil pour les surfaces
    $K3$}, Invent. Math. {\bf 15} (1972), 206-226

\bibitem[De79]{De2} P. Deligne,
\emph{Vari{\'e}t{\'e}s de Shimura: interpr{\'e}tation modulaire et
techniques de construction de mod{\`e}les canoniques},  in
{\it Automorphic Forms, Representations, and $L$-functions}
part. {\bf 2}, Proc. of Symp. in
Pure Math. {\bf 33}, American Mathematical Society (1979)  
247-290.

\bibitem[DMOS82]{Del82} P. Deligne, J. Milne, A. Ogus, K-Y. Shih, \emph{
    Hodge cycles, motives, and Shimura varieties}, LNM 900 (1982)

   
 
\bibitem[vdD98]{VDD} L. van den Dries, \emph{Tame Topology and
    o-minimal structures.} LMS lecture note series, \textbf{248},
Cambridge University Press, 1998.

\bibitem[vdDM94]{vdDM} L. van den Dries, C. Miller \emph{On the real
    exponential field with restricted analytic functions}, 
Israel J. Math. \textbf{85} (1994), 19--56.

\bibitem[FL81]{FL} E. Fortuna, S. \L ojasiewicz, \emph{Sur l'alg\'ebricit\'e des
    ensembles analytiques complexes}, J. Reine Angew. Math. {\bf 329}
  (1981) 215-220

\bibitem[Ga68]{Gabrielov} A.M Gabrielov, \emph{Projections of
    semi-analytic sets}, Funkt. Ana. i Prilozen \textbf{2} (1968), 18-30

  \bibitem[Gao16]{Gao} Z. Gao, \emph{Towards the Andr\'e-Oort
    conjecture for mixed Shimura varieties: the Ax-Lindemann theorem
    and lower bounds for Galois orbits of special points},
  J. Reine. Angew. Math. \textbf{732} (2017), 85-146

 \bibitem[Gao18]{Gao18} Z. Gao, \emph{Mixed Ax-Schanuel for the
     universal abelian varieties and some applications},
   Compos. Math. \textbf{156} (2020), no. 11, 2263-2297

   \bibitem[GK21]{GK} Z. Gao, B. Klingler, \emph{The Ax-Schanuel
       conjecture for variations of mixed Hodge structures}, \url{http://arxiv.org/abs/2101.10938}


  \bibitem[GGK12]{GGK} M. Green, P. Griffiths, M. Kerr, \emph{Mumford-Tate groups
    and domains. Their geometry and arithmetic},  Annals of
  Mathematics Studies \textbf{183}, Princeton University Press, 2012

\bibitem[G68]{Gri68} P. Griffiths, \emph{Period of integrals on algebraic
    manifolds, I}, Amer. J. Math., \textbf{90}, 568-626
  
\bibitem[G70a]{Grif}  P. Griffiths, \emph{Period of integrals on algebraic
    manifolds, III}, Inst. Hautes Etudes Sci. Publ. Math. \textbf{38}
  (1970) 125-180

  \bibitem[G70b]{Grifb} P. Griffiths, \emph{Periods of integrals on
      algebraic manifolds: Summary of main results and discussion of
      open problems}, Bull. Amer. Math. Soc. \textbf{76}, 228--296

\bibitem[GRT14]{GRT} P. Griffiths,  C. Robles, D. Toledo, {\em
    Quotients of non-classical flag domains are not algebraic},
  Algebr. Geom. \textbf{1} (2014), no. 1, 1-13

  \bibitem[GS69]{gsc} Griffiths P., Schmid W., Locally homogeneous complex
  manifolds, {\em Acta Math.} {\bf 123} (1969) 253-302

\bibitem[Gro84]{Gro} A. Grothendieck, \emph{Esquisse d'un programme} in 
Geometric Galois Actions vol. I, LMS Lecture Notes 242

\bibitem[GuVa04]{GuVa} S. Gukov, C. Vafa, \emph{Rational conformal field
    theories and complex multiplication}, Comm. Math. Phys. \textbf{246}
  (2004) no. 1, 181-210
  
\bibitem[HP12]{HP12} P. Habegger, J. Pila, \emph{Some unlikely
    intersections beyond Andr\'e-Oort}, Compos. Math. \textbf{148} (2012),
  no.1, 1-27
  
\bibitem[HP16]{HP16} P. Habegger, J. Pila, \emph{O-minimality and
    certain atypical intersections}, Ann. Sci. Ec. Norm. Super. (4),
  \textbf{49} (2016), no. 4, 813-858

\bibitem[Panorama]{panor} P. Habegger, G. R\'emond, T. Scanlon,
  E. Ullmo, A. Yafaev, \emph{Around the Zilber-Pink conjecture},
  Panoramas et Synth\`eses \textbf{52}, Soci\'et\'e Math\'ematique de
  France (2017)  

\bibitem[H51]{H51} W. Hodge, \emph{Differential forms on a K\"ahler
    manifold}, Proc. Cambridge Philos. Soc., \textbf{47}, (1951), 504-517

  
\bibitem[Ka85]{Ka} M. Kashiwara, \emph{The asymptotic behavior of a
    variation of polarized Hodge structure},
  Publ. Res. Inst. Math. Sci. \textbf{21} (1985), no. 4, 853-875

    

   
\bibitem[Khov80]{Khov} A. Khovanskii, \emph{On a class of systems of
    transcendental equations}, 
  Soviet Dokl. Math. \textbf{22} (1980) 762--765

  \bibitem[K17]{K17} B. Klingler, \emph{Hodge loci and atypical
    intersections: conjectures},
  \url{http://arxiv.org/abs/1711.09387}, accepted for publication in
  the book {\em Motives and complex multiplication}, Birkha\"user 

\bibitem[KO21]{KO} B. Klingler, A. Otwinowska, \emph{On the closure of
    the positive dimensional Hodge locus}, Inventiones Math. \textbf{225} (2021), no. 3, 857-883

  \bibitem[KOU]{KOU} B. Klingler, A. Otwinowska, D. Urbanik, \emph{On
      the fields of definition of Hodge loci},
    \url{http://arxiv.org/abs/2010.03359}, to appear in
    Ann. \'{E}c. Norm. Sup.
    
\bibitem[KUY16]{KUY} B. Klingler, E. Ullmo, A. Yafaev, \emph{The
    hyperbolic Ax-Lindemann-Weierstra\ss\- conjecture},
  Publ. Math. IHES \textbf{123},
  Issue 1, 333-360 (2016)

  \bibitem[KUY18]{KUY2}  B. Klingler, E. Ullmo, A. Yafaev, \emph{
    Bi-algebraic geometry and the Andr\'e-Oort conjecture: a survey},
  in Proceedings of 2015 AMS Summer
  Institute in Algebraic Geometry,  PSPMS 97-2, AMS, 2018, 319-360

  \bibitem[KP12]{KP12} J. K\'ollar, J. Pardon, \emph{Algebraic
      varieties with semialgebraic universal cover},
    J. Topol. \textbf{5} no.1, 199--212

\bibitem[Kr21a]{Kreutz1} T. Kreutz, \emph{Absolutely special
    subvarieties and absolute Hodge cycles},
  \url{http://arxiv.org/abs/2111.00216}

\bibitem[Kr21b]{Kreutz2} T. Kreutz, \emph{$\ell$-Galois special
    subvarieties and the Mumford-Tate conjecture}, \url{http://arxiv.org/abs/2111.01126}
  

\bibitem[Ma65]{Ma65} H.B. Mann, \emph{On linear relations between roots
        of unity}, Mathematik \textbf{12} (1965)
      
  \bibitem[MPT19]{MPT19} N.Mok, J. Pila, J. Tsimerman, \emph{Ax-Schanuel
      for Shimura varieties}, Ann. of Math. (2) \textbf{189} (2019),
    no.3, 945-978
    
    \bibitem[Moo98]{MoMo} B. Moonen, \emph{Linearity properties of Shimura
    varieties. I}, J. Algebraic Geom. \textbf{7} (1998), 539-567.


 \bibitem[Moore98]{Moore} G. Moore, \emph{Arithmetic and Attractors}, \url{https://arxiv.org/abs/hep-th/9807087}
  
 \bibitem[O18]{Orr} M. Orr, \emph{Height bounds and the Siegel
      property}, Algebra Number Theory {\bf 12} (2018), no. 2, 455-478.

    \bibitem[Oort94]{Oort} F. Oort, \emph{Canonical liftings and dense sets
    of CM-points}, Arithmetic Geometry (Cortona, 1994), 228-234,
  Symp. Math. XXXVII, Cambridge Univ. Press, Cambridge 1997
  
\bibitem[PS09]{PS} Y. Peterzil, S. Starchenko, \emph{Complex
    analytic geometry and analytic-geometric categories}, J. reine
  angew. Math. 626 (2009), 39-74
  
\bibitem[PS10]{PS1} Y. Peterzil, S. Starchenko \emph{Tame complex analysis and o-minimality.}
Proceedings of the ICM, Hyderabad, 2010. Available on first author's
web-page.

 \bibitem[PS13]{PS2} Y. Peterzil, S. Starchenko, \emph{
    Definability of restricted theta functions and families of abelian
    varieties}, Duke Math. J. \textbf{162}, (2013), 731-765


\bibitem[Pil11]{Pil} J. Pila, \emph{O-minimality and the Andre-Oort conjecture
    for $\CC^n$}, Ann. of Math. \textbf{173} (2011), 1779-1840.

\bibitem[Pil14]{Pil14} J. Pila, \emph{O-minimality and diophantine
    geometry}, Proceedings ot the ICM - Seoul (2014)

\bibitem[Pil15]{Pil15} J.Pila, \emph{Functional transcendence via
    o-minimality}, O-minimality and diophantine geometry, 66--99,
  London Math. Soc. Lecture Note Ser. \textbf{421}, Cambridge
  Univ. Press, Cambridge (2015)
  
\bibitem[PT14]{PT} J. Pila, J. Tsimerman, \emph{Ax-Lindemann for $\cA_g$.}
  Annals of Math \textbf{179} (2014), 659-681

\bibitem[PST21]{PST} J.Pila, A. Shankar, J. Tsimerman, with an
  appendix by H. Esnault and M. Groechenig, \emph{Canonical Heights on
    Shimura Varieties and the Andr\'e-Oort Conjecture}, \url{http://arxiv.org/abs/2109.08788}

\bibitem[PW06]{PW} J. Pila, A. Wilkie, \emph{The rational points on a
    definable set}, Duke Math. Journal \textbf{133}, (2006) 591-616.

\bibitem[PiZa08]{PZ} J. Pila, U. Zannier,
\emph{Rational points in periodic analytic sets and the Manin-Mumford
  conjecture}. Atti Accad. Naz. Lincei
Cl. Sci. Fis. Mat. Natur. Rend. Lincei (9) Mat. Appl. \textbf{19} (2008),
no. 2, 149-162.


\bibitem[Pink05]{Pink05} R. Pink, \emph{A combination of the
    conjectures of Mordell-Lang and Andr\'e-Oort}, in Geometric
 Methods in Algebra and Number Theory, Progress in Math. \textbf{253}
 Birkha\"user (2005), 251-282
 
\bibitem[Ray88]{Ray} M. Raynaud, \emph{Sous-vari\'et\'es d'une
    vari\'et\'e ab\'elienne et points de torsion}, in Arithmetic and
  Geometry, Vol. I, Progress in Math. \textbf{35} (1988)

  \bibitem[SaSc16]{SS} M. Saito, C. Schnell, \emph{Fields of definition of
        Hodge loci}, in Recent advances in Hodge theory, LMS \textbf{
        427} (2016) 275-291
  
\bibitem[Sc73]{Schmid} W. Schmid, \emph{Variation of Hodge structure:
    the singularities of the period mapping}, Invent. Math. \textbf{22}
  (1973), 211-319
  

  \bibitem[Se54]{Serre} J.P. Serre, \emph{G\'eom\'etrie alg\'ebrique
      et g\'eom\'etrie analytique}, Ann. Inst. Fourier (Grenoble)
    \textbf{6} (1955-56), 1--42

    
  \bibitem[ShWo95]{SW} H. Shiga, J. Wolfart, \emph{Criteria for complex
    multiplication and transcendence properties of automorphic
    functions}, J. Reine Angew. Math. \textbf{463} (1995) 1-25

  \bibitem[Ta20]{Tayou} S. Tayou, \emph{On the equidistribution of
      some Hodge loci}, J. Reine Angew. Math. \textbf{762} (2020),
    167-194

    \bibitem[TaTh21]{TayouTholozan} S. Tayou, N.Tholozan,
      \emph{Equidistribution of Hodge loci II}, \url{http://arxiv.org/abs/2103.15717}

  
  \bibitem[Tsi18a]{Tsi} J. Tsimerman, \emph{A proof of the Andr\'e-Oort
    conjecture for $\AAA_g$},  Ann. of Math. (2) \textbf{187} (2018),
  no. 2, 379-390

  \bibitem[Tsi18b]{Tsib} J. Tsimerman, \emph{Functional transcendence
      and arithmetic applications}, in Proc. Int. Cong. of
    Math. -2018, vol. 2, 453--472

\bibitem[U14]{U2} E. Ullmo, \emph{Applications du th\'eor\`eme
    d'Ax-Lindemann hyperbolique}, Compositio Math. \textbf{150} (2014),
  175-190 

\bibitem[U16]{Ull16} E. Ullmo, \emph{Structures sp\'eciales et
    probl\`eme de Zilber-Pink},  in ``Around the Zilber-Pink
  conjecture/Autour de la conjecture de Zilber-Pink'',
  Panor. Synth\`eses, \textbf{52}, Soc. Math. France, Paris (2017), 1--30

\bibitem[UY11]{UY1} E. Ullmo, A. Yafaev, \emph{A characterization of special subvarieties.} 
Mathematika \textbf{57} (2011) 263-273 

\bibitem[UY14b]{UY2} E. Ullmo, A. Yafaev, \emph{Hyperbolic Ax-Lindemann
    theorem in the cocompact case,} Duke Math. J. \textbf{163} (2014)
  433-463

  \bibitem[Ur21]{Urbanik} D. Urbanik, \emph{Absolute Hodge and $\ell$-adic
      monodromy}, \url{http://arxiv.org/abs/2011.10703}
    
\bibitem[V02]{Voisin02} C. Voisin, \emph{Hodge theory and complex
    algebraic geometry I}, Cambridge studies in advanced mathematics
  \textbf{76}, Cambridge University Press (2002)


  \bibitem[V07]{V07} C. Voisin, \emph{Hodge loci and absolute Hodge
      classes}, Compositio Math. \textbf{143}, 945-958, (2007)

\bibitem[Weil79]{Weil} A. Weil, \emph{Abelian varieties and the Hodge
ring}, Collected papers III, Springer Verlag, 421-429 (1979)
    
\bibitem[Wil96]{Wil} J. Wilkie, \emph{Model completeness results for
    expansions of the ordered field or real numbers by restricted
    Pfaffian functions and the exponential function},
  J. Amer. Math. Soc. \textbf{9} (1996), no.4, 1051-1094




\bibitem[YuZh]{YuZh} X. Yuan, S. Zhang, \emph{On the averaged Colmez
    conjecture}, Ann. of Math (2) \textbf{187} (2018), no.2, 533-638

\bibitem[Za12]{Zannier} U. Zannier, \emph{Some problems of unlikely
    intersections in arithmetic and geometry}, with appendices by
  David Masser, Annals of Mathematics Studies \textbf{181}, Princeton
  University Press (2012)
  
\bibitem[Zil02]{Zil02} B. Zilber, \emph{Exponential sums equations and
    the Schanuel conjecture}, J. London Math. Soc. (2) \textbf{65}
  (2002), no.1., 27--44
  
\end{thebibliography}
\end{document}